\documentclass[11pt]{article}
\usepackage{amssymb,amsmath}
\usepackage{color}
\usepackage[OT2, T1]{fontenc} 
\usepackage[utf8]{inputenc} 

\usepackage{latexsym,esint,
,mathrsfs} 
\usepackage[noadjust]{cite} 
\usepackage{verbatim,wasysym} 

\usepackage{graphicx} 

\usepackage[colorlinks=true, urlcolor=green!45!black, citecolor=red, linkcolor=blue,
pdfpagelabels, bookmarksnumbered, bookmarksopen]{hyperref} 

\usepackage[hyperpageref]{backref} 

\newcommand{\nequiv}{\not \equiv}

\usepackage[italian, english]{babel} 
\usepackage[OT2, T1]{fontenc} 

\usepackage{enumitem}

\def\claim#1.{\noindent {\bf #1.}}

\def\flushright#1{{\unskip\nobreak\hfil\penalty50\hskip2em\hbox{}\nobreak\hfil%
#1\parfillskip=0pt\finalhyphendemerits=0\par}}

\def\bull{\vrule height 1.8ex width 1.0ex depth .1ex }
\def\QED{\ifmmode\eqno\hbox{$\bull$}\else\flushright{\hbox{$\bull$}}\fi}

\newcommand{\parag}[1]{\left\{ \begin{aligned} #1 \end{aligned}\right.} 

\usepackage{scalerel}
\usepackage{tikz}
\usetikzlibrary{svg.path}
\definecolor{orcidlogocol}{HTML}{A6CE39}
\tikzset{
 orcidlogo/.pic={
 \fill[orcidlogocol] svg{M256,128c0,70.7-57.3,128-128,128C57.3,256,0,198.7,0,128C0,57.3,57.3,0,128,0C198.7,0,256,57.3,256,128z};
 \fill[white] svg{M86.3,186.2H70.9V79.1h15.4v48.4V186.2z}
 svg{M108.9,79.1h41.6c39.6,0,57,28.3,57,53.6c0,27.5-21.5,53.6-56.8,53.6h-41.8V79.1z M124.3,172.4h24.5c34.9,0,42.9-26.5,42.9-39.7c0-21.5-13.7-39.7-43.7-39.7h-23.7V172.4z}
 svg{M88.7,56.8c0,5.5-4.5,10.1-10.1,10.1c-5.6,0-10.1-4.6-10.1-10.1c0-5.6,4.5-10.1,10.1-10.1C84.2,46.7,88.7,51.3,88.7,56.8z};
 }
}
\newcommand\orcidicon[1]{\href{https://orcid.org/#1}{\mbox{\scalerel*{
\begin{tikzpicture}[yscale=-1,transform shape]
\pic{orcidlogo};
\end{tikzpicture}
}{|}}}}

\newtheorem{Theorem}{Theorem}[section]
\newtheorem{Proposition}[Theorem]{Proposition}
\newtheorem{Corollary}[Theorem]{Corollary}
\newtheorem{Lemma}[Theorem]{Lemma}
\newtheorem{Remark}[Theorem]{Remark}
\newtheorem{Definition}[Theorem]{Definition}

\newcommand{\R}{\mathbb{R}}
\newcommand{\N}{\mathbb{N}}
\newcommand{\Z}{\mathbb{Z}}

\newcommand{\mc}[1]{\mathcal{#1}}

\def\abs#1{|#1|}
\def\pabs#1{\left|{#1}\right|}
\def\norm#1{\|#1\|}
\def\eps{\varepsilon}
\newcommand{\dist}{{\rm dist}}

\begin{document}

\title
{Asymptotic decay of solutions for sublinear \\ fractional Choquard equations}

\author{
		\\Marco Gallo 
		\orcidicon{0000-0002-3141-9598}
		\\ \normalsize{Dipartimento di Matematica e Fisica}
		\\ \normalsize{Universit\`{a} Cattolica del Sacro Cuore}
		\\ \normalsize{Via della Garzetta 48, 25133 Brescia, Italy}
		\\ \normalsize{\href{mailto:marco.gallo1@unicatt.it}{\textcolor{black}{marco.gallo1@unicatt.it}}}
		\\
	}

\date{}

\maketitle


\begin{abstract} 
Goal of this paper is to study the asymptotic behaviour of the solutions of the following doubly nonlocal equation 
$$(-\Delta)^s u + \mu u = (I_{\alpha}*F(u))f(u) \quad \hbox{on $\mathbb{R}^N$}$$
where $s \in (0,1)$, $N\geq 2$, $\alpha \in (0,N)$, $\mu>0$, $I_{\alpha}$ denotes the Riesz potential and $F(t) = \int_0^t f(\tau) d \tau$ is a general nonlinearity with a sublinear growth in the origin. 
The found decay is of polynomial type, with a rate possibly slower than $\sim\frac{1}{|x|^{N+2s}}$.
The result is new even for homogeneous functions $f(u)=|u|^{r-2}u$, $r\in [\frac{N+\alpha}{N},2)$, and it complements the decays obtained in the linear and superlinear cases in \cite{DSS1, CGT3}. 
Differently from the local case $s=1$ in \cite{MS0}, new phenomena arise connected to a new ``$s$-sublinear'' threshold that we detect on the growth of $f$. 
To gain the result we in particular prove a Chain Rule type inequality in the fractional setting, suitable for concave powers.
\end{abstract}


\medskip


{
  \hypersetup{linkcolor=black} 
  \tableofcontents
}




\section{Introduction} 

The present paper is devoted to the study of the following doubly nonlocal equation
\begin{equation}\label{eq_main_intr}
(-\Delta)^s u + \mu u = \big(I_{\alpha}*F(u)\big)f(u) \quad \hbox{on $\R^N$}
\end{equation}
where $s \in (0,1)$, $N\geq 2$, $\alpha \in (0,N)$, $\mu>0$, $I_{\alpha}(x)=\frac{C_{N,\alpha}}{|x|^{N-\alpha}}$ is the Riesz potential and $(-\Delta)^s$ denotes the fractional Laplacian. 
The nonlinearity $F(t)=\int_0^t f(\tau) d\tau$ is assumed to be quite general, in the spirit of the papers by Berestycki and Lions \cite{BL1} and Moroz and Van Schaftingen \cite{MS2}, but the result is new even for the power case.
 In particular, we aim to study the asymptotic behaviour at infinity of the solutions: qualitative properties of this type have been already discussed when $f$ is linear or superlinear in \cite{CGT3} by the author, Cingolani and Tanaka, that is why we will restrict to the case of $f$ sublinear (in the origin).

Physically,
this doubly nonlocal model has different applications, in particular in the study of exotic stars: minimization properties related to \eqref{eq_main_intr} play indeed a fundamental role in the mathematical description of the dynamics of pseudo-relativistic boson stars \cite{ElSc0} and their gravitational collapse \cite{FJL}, as well as the evolution of attractive fermionic systems, such as white dwarf stars \cite{HLLS}. 
In fact, the study of the ground states to \eqref{eq_main_intr} gives information on the size of the critical initial conditions for the solutions of the corresponding pseudo-relativistic equation: 
in particular, when $s=\frac{1}{2}$, $N=3$, $\alpha=2$ and $f$ is a power, we obtain
$$\sqrt{-\Delta}u + \mu u = \left(\frac{1}{4 \pi r|x|}*|u|^r\right) |u|^{r-2} u \quad \hbox{in $\R^3$}$$
related to the so called 
\emph{massless boson stars equation} \cite{HL0}, where the pseudo-relativistic operator $\sqrt{-\Delta + m^2}-m$ collapses to the square root of the Laplacian; we refer to \cite{GalT} and references therein for a soft introduction. 
Other applications can be found in relativistic physics and in quantum chemistry \cite{ArMe,DOS
} and in the study of graphene \cite{LMM}. 

\medskip

Mathematically, when $s=1$ and $f$ is a power, 
that is
\begin{equation}\label{eq_MS}
-\Delta u + \mu u = \big(I_{\alpha}*|u|^r\big)|u|^{r-2} u \quad \hbox{in $\R^N$},
\end{equation}
Cingolani, Clapp and Secchi in \cite[Proposition A.2]{CCS1} obtained an exponential decay of positive solutions whenever $r\geq 2$, which means that the effect of the classical Laplacian prevails.
 Afterwards Moroz and Van Schaftingen in \cite{MS0} (see also \cite{
 MS3} and \cite{CLO, ClSa}) extended the previous analysis in the case of ground state solutions to all the possible values of $r$ in the range $[\frac{N+\alpha}{N}, \frac{N+\alpha}{N-2}]$, in particular by finding a polynomial decay when $f$ is sublinear (i.e., the Choquard term effect prevails).
%
They prove the following result \cite[Theorem 4]{MS0}.

\begin{Theorem}[\cite{MS0}]\label{thm_MVS}
Let $u \in H^1(\R^N)$ be a nonnegative ground state of \eqref{eq_MS}, and $r \in [\frac{N+\alpha}{N}, \frac{N+\alpha}{N-2}]$.
Assume $\mu=1$. 
Then
\begin{itemize}
\item if $r>2$, then
$$\lim_{|x|\to +\infty} u(x) |x|^{\frac{N-1}{2}} e^{|x|} \in (0, +\infty);$$
\item if $r=2$, then 
$$\lim_{|x|\to +\infty} u(x) |x|^{\frac{N-1}{2}} e^{\int_{\nu}^{|x|} \sqrt{1- \frac{\nu^{N-\alpha}}{t^{N-\alpha}}} dt } \in (0, +\infty)$$
for some explicit $\nu=\nu(u)$;
\item if $r<2$, then
$$ \lim_{|x|\to +\infty} u(x) |x|^{\frac{N-\alpha}{2-r}} = C(N, \alpha, r, u) \in (0, +\infty)$$
where 
\begin{equation}\label{eq_costant_stima}
C(N,\alpha, r, u) := \big( C_{N,\alpha} \norm{u}_r^r\big)^{\frac{1}{2-r}}
\end{equation}
with $C_{N,\alpha}:=\frac{\Gamma(\frac{N-\alpha}{2})}{2^{\alpha} \pi^{N/2} \Gamma(\frac{\alpha}{2})}$.
\end{itemize}
\end{Theorem}

Notice that, when $\mu\neq 1$, the frequency $\mu$ influences both the limiting constants and -- when $r \geq 2$ -- the speed of the exponential decays. 
We refer also to \cite[Section 8.2]{DG0} for some results on convolution equations with non-variational structure.

\smallskip

The case of the fractional Choquard equation $s \in (0,1)$ with homogeneous $f$, that is
\begin{equation}\label{eq_fr_ch_hom}
(-\Delta)^s u + \mu u = \big(I_{\alpha}*|u|^r\big)|u|^{r-2} u \quad \hbox{in $\R^N$},
\end{equation}
 has been studied by D'Avenia, Siciliano and Squassina in \cite{DSS1} (see also 
 \cite{MZ0, BBS, ZY0}
for other related results). 
In this paper the authors gain existence of ground states, multiplicity and qualitative properties of solutions: 
in particular they obtain asymptotic decay of solutions whenever the source is linear or superlinear, that is when $r\geq 2$ (see also \cite{BBMP} for the $p$-fractional Laplacian counterpart); in this case the rate is polynomial, as one can expect dealing with the fractional Laplacian. More specifically, it does not depend on $\alpha$, and they prove the following theorem.

\begin{Theorem}[\cite{DSS1}]\label{thm_DSS1}
Let $u\in H^s(\R^N)$ be a solution of \eqref{eq_fr_ch_hom}, and assume $r \in [2, \frac{N+\alpha}{N-2s}]$. Then
\begin{equation}\label{eq_stima_superlin}
0 < \liminf_{|x|\to +\infty} |u(x)| |x|^{N+2s} \leq \limsup_{|x|\to +\infty} |u(x)| |x|^{N+2s} <+\infty.
\end{equation}
\end{Theorem}


In this paper, we aim to study the fractional Choquard case $s\in (0,1)$, $\alpha \in (0,N)$, in presence of general, sublinear nonlinearities. 
We point out that the arguments in \cite{MS0} cannot be directly adapted to the fractional framework: for instance, we see that the explicit computation of the fractional Laplacian of some comparison function is not possible, and the choice of the comparison functions itself is hindered by some growth condition typical of the fractional framework; moreover, it is not obvious that all the weak solutions are pointwise solutions, and neither one can deduce that the concave power of a pointwise solution is indeed a solution (of a different equation) itself.

\smallskip

We start by presenting the case of homogeneous powers $f$, which has an interest on its own.
Since in the superlinear case the rate of convergence is of the type $\sim \frac{1}{|x|^{N+2s}}$, in the sublinear case we generally expect a slower decay. 
Actually this is what we find, as the following theorem states. 
\begin{Theorem}\label{thm_homog0_ws}
Let $u\in H^s(\R^N)$, strictly positive, radially symmetric and decreasing, be a weak solution of \eqref{eq_fr_ch_hom}.
 Let $r \in [\frac{N+\alpha}{N}, 2)$ and set 
 \begin{equation}\label{eq_defin_beta}
 \beta:= \min\left\{ \frac{N-\alpha}{2-r}, N+2s\right\} \in [N, N+2s]. 
 \end{equation}
Then
$$0 < \liminf_{|x| \to +\infty} u(x)|x|^{\beta} \leq \limsup_{|x|\to +\infty} u(x) |x|^{\beta} < + \infty .$$
Moreover, in the case $r \in [\frac{N+\alpha}{N}, \frac{N+\alpha+4s}{N+2s})$ (i.e. $\beta < N+2s$), 
we have the \emph{sharp decay}
\begin{equation}\label{eq_corol_sharp_decay}
 \lim_{|x|\to +\infty} u(x) |x|^{\frac{N-\alpha}{2-r}} = \left( \frac{C_{N,\alpha}\norm{u}_r^r}{\mu} \right)^{\frac{1}{2-r}}.
 \end{equation}
\end{Theorem}
We notice that, if $\mu=1$, the constant in \eqref{eq_corol_sharp_decay} is coherent with \eqref{eq_costant_stima}.
We refer to Remark \ref{rem_more_gen} 
for some comments and generalizations on the assumptions. 
 This result in particular applies to ground states solutions (see Definition \ref{def_pohozaev_min2}). 

\begin{Corollary}\label{corol_homog0_gs}
Let $u$ be a positive ground state of \eqref{eq_fr_ch_hom}. 
Then the conclusions of Theorem \ref{thm_homog0_ws} hold.
\end{Corollary}

We highlight that the found decay of the ground states might give information, when $r<2$, also on the twice Gateaux differentiability of the corresponding functional and on the nondegenaracy of the ground state solution itself, see \cite{MS0} (see also \cite[Section 3.3.5]{MS3}). 
Moreover 
this information on the decay may be exploited to study fractional Choquard equations with potentials $V=V(x)$ approaching, as $|x|\to +\infty$, some $V_{\infty}>0$ from above or oscillating, in the spirit of \cite{
MPS2}. It might be further used, for example, in the semiclassical analysis of concentration phenomena, see e.g. \cite{CG0}.

\smallskip

In both the estimates from above and below in Theorem \ref{thm_homog0_ws} we rely on some comparison principle and the use of some auxiliary function whose fractional Laplacian is related to the Gauss hypergeometric function.
For the estimate from above we succeed in working with the weak formulation of the problem; 
on the other hand, in order to deal with the estimate from below, we find the necessity of working with $u^{2-r}$, where $2-r \in (0,1)$: this concave power of the solution may fail to lie in $H^s(\R^N)$, and thus we cannot treat the problem with its weak formulation. 
The pointwise formulation seems to arise some problems as well, since the fractional Laplacian of $u^{2-r}$ needs some restrictive assumption on $\alpha,s,N$ and $r$ in order to be well defined. 
This is why we work with a viscosity formulation of the problem, obtaining a C\'ordoba-C\'ordoba type inequality for concave functions (see Lemma \ref{lem_chain_rule}). 
We remark that the estimate from above may be treated with the viscosity formulation as well.

\medskip

The paper is organized as follows.
In Section \ref{sec_extra_referee} we make some comments on the found results and present some generalizations, in particular for the case of a general nonlinearity $f=f(t)$ in \eqref{eq_main_intr}.
In Section \ref{sec_prelim} we introduce definitions and notations, collecting some existence and comparison results in Appendix \ref{app_ex_comp}. In Section \ref{sec_frac_aux} we introduce some suitable auxiliary function (see Appendix \ref{app_frac_dec} for some related asymptotic property) and establish some asymptotic behaviour on suitable comparison functions; other preliminary estimates are studied in Section \ref{sec_prel_est}.
Then in Section \ref{sec_above} we deal with the estimate from above, by working with the weak formulation, while in Section \ref{sec_est_bel_vis} we study the asymptotic behaviour from below, by exploiting a viscosity formulation and proving a fractional Chain Rule, suitable for concave functions. 
Finally in Section \ref{sec_proof_main_2} we conclude the proofs of the main results.

\section{Comments and generalizations}
\label{sec_extra_referee}

Joining the results in Theorem \ref{thm_DSS1} and Theorem \ref{thm_homog0_ws} we obtain the following picture of the asymptotic decay of fractional Choquard equations.
\begin{Corollary}
Let $u$ be a positive ground state of \eqref{eq_fr_ch_hom}, with $r \in [ \frac{N+\alpha}{N}, \frac{N+\alpha}{N-2s}]$. 
\begin{itemize}
\item If $r \in [\frac{N+\alpha+4s}{N+2s}, \frac{N+\alpha}{N-2s}]$, then
$$0 < \liminf_{|x|\to +\infty} u(x) |x|^{N+2s} \leq \limsup_{|x|\to +\infty} u(x) |x|^{N+2s} <+\infty.$$
\item If $r \in [\frac{N+\alpha}{N}, \frac{N+\alpha+4s}{N+2s}]$, then 
$$0 < \liminf_{|x|\to +\infty} u(x) |x|^{\frac{N-\alpha}{2-r}} \leq \limsup_{|x|\to +\infty} u(x) |x|^{\frac{N-\alpha}{2-r}} <+\infty;$$
in particular, $\frac{N-\alpha}{2-r}=N$ in the lower critical case $r= \frac{N+\alpha}{N}$.
\end{itemize}
\end{Corollary}
By the previous Corollary we see that the exponent 
$$r^*_{s,\alpha}:= \frac{N+\alpha+4s}{N+2s},$$
$r^*_{s,\alpha} \in (\frac{N+\alpha}{N}, 2)$, separates the cases where the fractional Laplacian influences more the rate of convergence (which does not depend on $\alpha$), from the cases where the asymptotic behaviour is dictated by the Choquard term (which does not depend on $s$). 
This phenomenon seems to highlight a difference between the fractional and the local case, where the separating exponent is $r=2$ (see Theorem \ref{thm_MVS}): indeed, when $r \in \left(r^*_{1,\alpha}, 2\right)$, the arbitrary big (as $r \to 2$) polynomial behaviour $\sim \frac{1}{|x|^{\frac{N-\alpha}{2-r}}}$ keeps being slower than the exponential decay induced by the classical Laplacian; this is not the case when compared with the polynomial decay induced by the fractional Laplacian, and this is why this new phenomenon appears in this range. 
Thus $r^*_{s,\alpha}$ can ben seen as a kind of \emph{$s$-subquadratic} threshold for the growth of $F$; 
set instead
$$p^*_{s,\alpha}:=r^*_{s,\alpha}-1 = \frac{\alpha+2s}{N+2s},$$
it can be seen as a \emph{$s$-sublinear} threshold for the growth of $f$. 
Notice that
$$ r^*_{s,\alpha} \stackrel{s \to 0} \to \frac{N+\alpha}{N}, \quad r^*_{s,\alpha} \stackrel{\alpha \to N} \to 2,$$
while
$$ r^*_{s,\alpha} \stackrel{s \to 1} \to \frac{N+\alpha+4}{N+2} \in \Big(\frac{N+\alpha}{N}, 2\Big), \quad r^*_{s,\alpha} \stackrel{\alpha \to 0} \to \frac{N+4s}{N+2s} \in (1,2).$$
%
It might be interesting to investigate other possible phenomena on fractional Choquard equations when $r$ is above and below this exponent $r^*_{s,\alpha}$, or also possible phenomena in $(r^*_{1,\alpha},2)$ for the local Choquard equation.
We refer also to the recent paper \cite[Theorem 1.4]{Gre0} where asymptotic decay results are studied in a different framework (still involving the fractional Laplacian and the Riesz potential); here a threshold different from the classical case $s=1$ is detected as well.

%

\begin{Remark}
We notice that, fixed a positive solution $u$, by setting
$$\rho:= I_{\alpha}* u^r$$
equation \eqref{eq_fr_ch_hom} can be rewritten as
$$(-\Delta)^s u + \mu u = \rho(x) u^{r-1}.$$
When $\mu=0$ and $\rho(x) \leq \frac{1}{|x|^{\gamma}}$ with $\gamma > N$, this fractional sublinear equation ($r\in (0,2)$) has been studied in \cite{PT0} (see also \cite[Theorem 4.4]{GM0} where they extend the result to $\gamma>2s$): here the authors find an estimate from above of the asymptotic decay of the solutions, which is strictly slower than $\sim \frac{1}{|x|^N}$. 
Notice that, in our case, $\rho= I_{\alpha}*u^r$ decays at most as $\sim \frac{1}{|x|^{N-\alpha}}$ (see \cite[Lemma 4.6]{GM0} and
\cite[page 801]{MS3}) 
and we discuss the strict positivity of $\mu$. See also \cite{DSS1, Le0} for more results on the zero mass case. 


\end{Remark}


We pass now to more general nonlinearities, and study \eqref{eq_main_intr}. For the whole paper we assume the following conditions on $f$ in order to give sense to appearing integrals: 
\begin{itemize}
\item[\textnormal{(f1)}] \label{(f1)}
 $f \in C(\R, \R)$, $F(t) = \int_0^t f(\tau) d \tau$;
\item[\textnormal{(f2)}] \label{(f2)}
$f$ satisfies
 $$i) \; \limsup_{t \to 0} \frac{|tf(t)|}{|t|^{\frac{N+ \alpha}{N}}} <+\infty, \quad
 ii) \; \limsup_{ |t| \to + \infty} \frac{|t f(t)|}{|t|^{\frac{N+ \alpha}{N-2s}}} <+\infty,$$
or equivalently there exists $C >0$ such that for every $t \in \R$, 
$$|t f(t)| \leq C \big(|t|^{\frac{N + \alpha}{N}} + |t|^{\frac{N+ \alpha}{N-2s}}\big).$$
\end{itemize}
In particular, \hyperref[(f2)]{\textnormal{(f2)}} implies
 \begin{equation}\label{eq_condiz_F}
i) \; \limsup_{t \to 0} \frac{|F(t)|}{|t|^{\frac{N+ \alpha}{N}}} <+\infty, \quad
 ii) \; \limsup_{ |t| \to + \infty} \frac{|F(t)|}{|t|^{\frac{N+ \alpha}{N-2s}}} <+\infty,
\end{equation}
or equivalently that there exists $C >0$ such that for every $t \in \R$, 
$$|F(t)| \leq C \big(|t|^{\frac{N + \alpha}{N}} + |t|^{\frac{N+ \alpha}{N-2s}}\big).$$
These conditions have been introduced in \cite{MS2} 
for the local case $s=1$, extending \cite{BL1} where the seminal case of local nonlinearities is treated. 
These critical exponents have then been adapted to the fractional case $s \in (0,1)$ in \cite{DSS1}, while the general case \hyperref[(f1)]{\textnormal{(f1)}}-\hyperref[(f2)]{\textnormal{(f2)}} has been introduced in \cite{CGT2}. 
This set of assumptions covers different types of nonlinearities, such as pure powers, both odd $f(u)=|u|^{r-1} u$ or even $f(u)=|u|^r$, combination of powers $f(u)= u^r \pm u^q$ (standing for cooperation or competition), asymptotically linear (saturable) nonlinearities $\frac{u^{r+1}}{1+u^r}$ (which appear in nonlinear optics \cite{DLWZPH}) and many others. 
Notice that these assumptions include the case of critical nonlinearities, both in the origin and at infinity. 

In the papers \cite{CGT2, CGT3, CG1,CGT5} (see also \cite{GalT}) the authors study existence and multiplicity of normalized solutions and of Pohozaev minima for \eqref{eq_main_intr}, as well as qualitative properties of solutions, such as regularity, positivity, radial symmetry and Pohozaev identities. 
In particular in \cite{CGT3} they extend Theorem \ref{thm_DSS1} to the case of general nonlinearities, by proving the polynomial asymptotic behaviour of solutions whenever $f$ is linear or superlinear \emph{in the origin}. 
That is, by assuming 
$\limsup_{t \to 0} \frac{|f(t)|}{|t|}<+\infty$ 
they gain that every positive weak solution $u$ satisfies \eqref{eq_stima_superlin}.

\smallskip

In this paper, we further investigate the asymptotic behaviour of the solutions of the fractional Choquard equation \eqref{eq_main_intr} when $f$ is \emph{sublinear in the origin}. 
Thus we consider the following additional assumptions:

\begin{itemize}
\item[\textnormal{(f3)}] \label{(f3)}
 there exists $r \in [\frac{N+\alpha}{N}, 2)$ such that
$$\limsup_{t \to 0^+} \frac{|f(t)|}{t^{r-1}} \in [0, +\infty),$$
i.e., for some $\bar{C}>0$ and $\delta \in (0,1)$ we have 
\begin{equation}\label{eq_cond_2sublin_f}
|f(t)| \leq \bar{C} t^{r-1} \quad \hbox{for $t\in (0,\delta)$};
\end{equation}
\item[\textnormal{(f4)}] \label{(f4)} 
there exists $r \in [\frac{N+\alpha}{N}, 2)$ such that
$$\liminf_{t \to 0^+} \frac{f(t)}{t^{r-1}} \in (0, +\infty),$$
i.e., for some $\underline{C}>0$ and $\delta \in (0,1)$ we have 
\begin{equation}\label{eq_cond_4sublin_f}
f(t) \geq \underline{C} t^{r-1} \quad \hbox{for $t \in (0,\delta)$}.
\end{equation}
\end{itemize}

\noindent
A sufficient condition for \hyperref[(f3)]{\textnormal{(f3)}} 
is clearly given by
\begin{equation}\label{eq_strongf3}
\limsup_{t \to 0^+} \frac{f(t)}{t^{r-1}}=0 \quad \hbox{for some $r \in [\frac{N+\alpha}{N}, 2)$},
\end{equation}
which means that $\bar{C}$ can be taken arbitrary small in \eqref{eq_cond_2sublin_f} (up to taking $\delta$ sufficiently small); in particular it includes logarithmic nonlinearities $f(t)=t \log(t^2)$, 
where $r$ can be chosen arbitrary close to $2$. 
A sufficient condition for \hyperref[(f4)]{\textnormal{(f4)}} is instead given (for example) by a local Ambrosetti-Rabinowitz condition ($f(t) t \geq r F(t) >0$ for $t \in (0,\delta)$).
The restriction in \hyperref[(f3)]{\textnormal{(f3)}} and \hyperref[(f4)]{\textnormal{(f4)}} to right neighborhoods of zero is due to the fact we deal with positive solutions.

\smallskip

We eventually come up with the following generalization of Theorem \ref{thm_homog0_ws}.
\begin{Theorem}\label{thm_main}
Assume \hyperref[(f1)]{\textnormal{(f1)}}-\hyperref[(f2)]{\textnormal{(f2)}}, and let $u\in H^s(\R^N)$, strictly positive, radially symmetric and decreasing, be a weak solution of \eqref{eq_main_intr}. Let $r \in [\frac{N+\alpha}{N}, 2)$ and $\beta$ as in \eqref{eq_defin_beta}.
\begin{itemize}
\item[\textnormal{(i)}] 
Assume \hyperref[(f3)]{\textnormal{(f3)}}. 
Then 
$\limsup_{|x|\to +\infty} u(x) |x|^{\beta} \in (0, +\infty).$
\item[\textnormal{(ii)}] 
Assume \hyperref[(f4)]{\textnormal{(f4)}}, $f$ locally H\"older continuous and $\int_{\R^N}F(u)>0$ (e.g. $F\geq 0$ on $(0,+\infty)$). 
Then 
$\liminf_{|x| \to +\infty} u(x)|x|^{\beta} \in (0, +\infty).$
\end{itemize}
If both conditions in (i) and (ii) hold, together with $\overline{C}=\underline{C}$ (i.e., $f$ is a power near the origin) and $r \in [\frac{N+\alpha}{N}, \frac{N+\alpha+4s}{N+2s})$, then 
we have the sharp decay
\begin{equation}\label{eq_sharp_decay_F} 
\lim_{|x|\to +\infty} u(x) |x|^{\frac{N-\alpha}{2-r}} = \left( \frac{C_{N,\alpha} \big(\lim_{t\to 0^+} \frac{f(t)}{t^{r-1}}\big) \int_{\R^N} F(u)}{\mu } \right)^{\frac{1}{2-r}} 
\end{equation}
where $C_{N,\alpha}>0$ is given in \eqref{eq_costant_stima}.
\end{Theorem}
\begin{Remark}\label{rem_more_gen}
We highlight that the conclusions of Theorem \ref{thm_main} (as well as of Theorem \ref{thm_homog0_ws}) hold in more general cases. 
\begin{itemize}
\item The case
$$\lim_{t\to 0^+} \frac{f(t)}{t}=+\infty$$
in a non-strict sense (i.e. $\lim_{t\to 0} \frac{|f(t)|}{|t|^{r-1}}=0$ for each $r \in [\frac{N+\alpha}{N}, 2)$, for example $f(t) \sim - t \log(t^2)$) is included, and as we expect the decay is of order $\sim \frac{1}{|x|^{N+2s}}$. It is sufficient to apply the argument of Remark \ref{rem_dec_2sN} (since $f(t) \geq \underline{C} t$ for $t$ small and positive), and the results in Proposition \ref{prop_estim_above} (after having chosen a whatever $r \in [r^*_{\alpha,s}, 
2)$). 
\item The conclusions hold also without assuming radial symmetry and monotonicity of $u$, but by assuming a priori that
$$\limsup_{|x|\to +\infty} |u(x)| |x|^{\omega} < +\infty$$
for some $\omega> \frac{N^2}{N+\alpha}$: see Remark \ref{rem_cond_asym_dec}. When $u\in L^q(\R^N)$, $q < \frac{N+\alpha}{N}$, is radially symmetric and decreasing, this is the case with $\omega=\frac{N}{q}$ (see Remark \ref{rem_dec_N}); in particular, if $q=1$, we have $\omega=N$. Notice that $u$ is automatically radially symmetric and decreasing when $u\in C^{1,1}_{loc}(\R^N)$, $f(u)=|u|^{r-2}u$ and $\omega>\frac{\alpha}{r-1}$ thanks to \cite[Theorem 1]{Le1} (see also \cite[Theorem 1.3]{WY0}). 

\item In light of the previous remark, we highlight that the estimate from above actually holds true also for nonnegative solutions $u \geq 0$; see Proposition \ref{prop_estim_above}; moreover, it can be further extended to $|u|$ in the case of changing sign solutions, by applying a Kato's inequality \cite[Theorem 3.2]{Amb0}.
\item The conclusions hold also for solutions $u\in L^1(\R^N) \cap C(\R^N)$ in the viscosity sense, without assuming $f$ H\"older continuous (which is needed in (ii) only to pass from weak to viscosity solutions): see Section \ref{sec_est_bel_vis}.
\item When \hyperref[(f4)]{\textnormal{(f4)}} holds, we actually have $F(t) \geq \underline{C} \frac{t^r}{r}$ for $t \in (0,\delta)$; thus, 
being also $u\in L^{\infty}(\R^N)$, the condition $\int_{\R^N} F(u)>0$ means that $F$ is not \emph{too negative} in $[\delta, \norm{u}_{\infty}]$. We highlight that the energy term $\int_{\R^N} \big(I_{\alpha}*F(u)\big) F(u)$ is always positive (see e.g. \cite{CGT5}).
\item We find some estimates on the asymptotic constants, which are coherent, when $r \in [\frac{N+\alpha}{N}, 
r^*_{s,\alpha})$, 
with the one found in Theorem \ref{thm_MVS} and Theorem \ref{thm_main}: 
see Propositions \ref{prop_estim_above} and \ref{prop_below_2}. 
We notice that \eqref{eq_fr_ch_hom} is obtained by \eqref{eq_main_intr} formally choosing $f(t)=\sqrt{r}|t|^{r-2}t$. In the paper -- up to well posedness and regularity -- we do not use that $F$ is the primitive of $f$: in particular, we do not apply \hyperref[(f3)]{\textnormal{(f3)}} and \hyperref[(f4)]{\textnormal{(f4)}} to $F$. Thus we can arbitrary move constants from $f$ to $F$ in our arguments to adjust -- for example -- the value of $\underline{C}$, and this allows to gain the result for every $\mu>0$. 
\end{itemize}
\end{Remark}

\smallskip

Our results apply in particular to Pohozaev minima of the equation (see Definition \ref{def_pohozaev_min2}), whenever some symmetric assumption is assumed on $f$, that is
\begin{itemize}
\item[\textnormal{(f5)}] \label{(f5)}
$f$ is odd or even, with constant sign on $(0,+\infty)$ and locally H\"older continuous.
\end{itemize}
We refer to \cite{CG1} for discussions on the assumption \hyperref[(f5)]{\textnormal{(f5)}}. We notice that, since every Pohozaev minimum has strict constant sign \cite{CG1}, it is not restrictive to assume a priori the sign of $u$.

\begin{Corollary}\label{corol_main}
Assume \hyperref[(f1)]{\textnormal{(f1)}}-\hyperref[(f2)]{\textnormal{(f2)}} and \hyperref[(f5)]{\textnormal{(f5)}}. Let $u$ be a (positive) Pohozaev minimum of \eqref{eq_main_intr}. Then the conclusions of Theorem \ref{thm_main} hold.
\end{Corollary}


We finally want to highlight that our results may be adapted to the local case $s=1$, extending Theorem \ref{thm_MVS} to general nonlinearities, studied in \cite{MS2}. 
We leave the details to the reader, observing that in this case the rate of decaying is simply given by $\beta= \frac{N-\alpha}{2-r}$, since, as already observed, the solutions of the homogeneous linear (associated) equation decay exponentially.

\begin{Theorem}\label{thm_main_loc}
Let $s=1$, and assume \hyperref[(f1)]{\textnormal{(f1)}}-\hyperref[(f2)]{\textnormal{(f2)}} (where the upper critical exponent is substituted by $\frac{N+\alpha}{N-2}$). Let $u\in H^1(\R^N)$, strictly positive, radially symmetric and decreasing, be a solution of 
$$-\Delta u + \mu u = \big(I_{\alpha}*F(u)\big)f(u) \quad \hbox{on $\R^N$};$$
in particular, $u$ may be a ground state. Let $r \in [\frac{N+\alpha}{N}, 2)$.
\begin{itemize}
\item[\textnormal{(i)}] Assume \hyperref[(f3)]{\textnormal{(f3)}}. 
Then 
$\limsup_{|x|\to +\infty} u(x) |x|^{ \frac{N-\alpha}{2-r}} \in (0, +\infty).$
\item[\textnormal{(ii)}] Assume \hyperref[(f4)]{\textnormal{(f4)}} and $\int_{\R^N}F(u)>0$. 
Then
$\liminf_{|x| \to +\infty} u(x)|x|^{\frac{N-\alpha}{2-r}} \in (0, +\infty).$
\end{itemize}
If both conditions (i) and (ii) hold, together with $\overline{C}=\underline{C}$, then \eqref{eq_sharp_decay_F} holds.
\end{Theorem}

\section{Preliminaries}

\subsection{Definitions and notations}
\label{sec_prelim}

Let $s \in (0,1)$ and $\alpha \in (0, N)$, where $N\geq 2$. 
We will denote by $C^{k,\sigma}(\R^N)$ the space of the functions in $C^k(\R^N)$ with $\sigma$-H\"olderian $k$-derivatives, and more briefly we will write $C^{\gamma}(\R^N):= C^{[\gamma], \gamma-[\gamma]}(\R^N)$ for any $\gamma>0$. 
The same notations apply to the local case $C^{\gamma}_{loc}(\R^N)$. Moreover we write $\norm{\cdot}_{p}=\norm{\cdot}_{L^p(\R^N)}$ for the classical $L^p$ norm in the entire space, $p \in [1, +\infty]$, and we will use also the following notation
$$\norm{f}_{\infty, \theta}:= \norm{f(\cdot )(1+|\cdot|^{\theta})}_{\infty}$$
for any $\theta>0$. 
Finally by $f \sim g$ as $|x|\to +\infty$ we mean that $\lim_{|x|\to +\infty} \frac{f(x)}{g(x)}=1$.

Let the fractional Laplacian be defined via Fourier transform \cite{DPV}
$$(-\Delta)^s u = \mc{F}^{-1}(|\xi|^{2s} \mc{F}(u)),$$
while, when $u$ is regular enough, we can write \cite[Proposition 3.3]{DPV}
$$(-\Delta)^s u(x) = C_{N,s} \int_{\R^N} \frac{u(x)-u(y)}{|x-y|^{N+2s}} dy, \quad x \in \R^N$$
where $C_{N,s}:=\frac{4^s \Gamma(\frac{N+2s}{2})}{\pi^{N/2} |\Gamma(-s)|}>0$ and the integral is in the principal value sense. 
A sufficient condition in order to have $(-\Delta)^s u$ well defined pointwise is given by \cite[Proposition 2.4]{Sil0} (see also \cite[Proposition 2.15]{Gar0} and \cite[Proposition 2.1]{CGT5}).
\begin{Proposition}[Pointwise well posedness]\label{prop_well_posed}
Let $x_0 \in \R^N$. Then, if
$u \in L^p(\R^N)\cap C^{\gamma}(U)$ for some $p \in [1,+\infty]$, $\gamma >2s$ and $U$ open neighborhood of $x_0$,
then $(-\Delta)^s u(x_0)$ is well defined. Moreover, 
$(-\Delta)^s u \in C(U)$. 
\end{Proposition}

We introduce, for any $\Omega \subset \R^N$ and $s \in (0,1)$,
$$H^s(\Omega) := \left\{ u \in L^2(\Omega) \mid [u]_{H^s(\Omega)}^2:=\int_{\Omega} \int_{\Omega} \frac{|u(x)-u(y)|^2}{|x-y|^{N+2s}}dy < +\infty \right \},$$
endowed with
$$\norm{u}_{H^s(\Omega)}^2:= \norm{u}_{L^2(\Omega)}^2 + [u]_{H^s(\Omega)}^2.$$
We recall that \cite[Theorem 5.4 and 6.7]{DPV}, when $\Omega$ is for example an open set with $C^{0,1}$ bounded boundary, we have $H^s(\Omega) \hookrightarrow L^2(\Omega) \cap L^{2^*_s}(\Omega)$, where $2^*_s:= \frac{2N}{N-2s}$. 
Moreover we set 
$$H^s_{loc}(\R^N):= \left\{ u: \R^N \to \R \mid u \in H^s(\Omega) \hbox{ for each $\Omega \subset \subset \R^N$}\right\}$$
and \cite[Section 4.3.2]{Tri0}
$$X^s_0(\Omega):= \left\{ w \in H^s(\R^N) \mid w=0 \; \hbox{on $\Omega^c$}\right\}.$$
In the case $\Omega = \R^N$ we also have the following relation \cite[Proposition 3.6]{DPV} 
$$\norm{(-\Delta)^{s/2} u}_{L^2(\R^N)}^2 = \tfrac{1}{2} C_{N,s} [u]_{H^s(\R^N)}^2$$
which leads to the following formulation via Fourier transform 
$$H^s(\R^N) = \left\{ u \in L^2(\R^N) \mid |\xi|^{s} \widehat{u} \in L^2(\R^N) \right \};$$ 
this definition extends also to every $s>0$ \cite{FQT}.

We further recall the Riesz potential 
\begin{equation}\label{eq_def_Riesz}
I_{\alpha}(x) = \frac{C_{N,\alpha}}{|x|^{N-\alpha}}, \quad x \in \R^N \setminus \{0\}
\end{equation}
where $C_{N,\alpha}:=\frac{\Gamma(\frac{N-\alpha}{2})}{2^{\alpha} \pi^{N/2} \Gamma(\frac{\alpha}{2})}>0$: by the Hardy-Littlewood-Sobolev inequality we have
$$f \in L^r(\R^N) \mapsto I_{\alpha}*f \in L^h(\R^N)$$
continuous whenever $r,h \in (1,+\infty)$ satifsy $\frac{1}{r}-\frac{1}{h}=\frac{\alpha}{N}$. 

\begin{Remark}\label{rem_conv_welldef}
Arguing as in \cite[Proposition 4.5]{CGT3} we see that $I_{\alpha}*F(u) \in C_0(\R^N)$ (and thus it is well defined pointwise) if $F(u)$ lies in $L^{\frac{N}{\alpha}-\eps}(\R^N) \cap L^{\frac{N}{\alpha}+\eps}(\R^N)$ for some $\eps>0$. 
%
%
Assuming \hyperref[(f1)]{\textnormal{(f1)}}-\hyperref[(f2)]{\textnormal{(f2)}} on $f$, we need to assume that $u \in L^{\frac{N+\alpha}{\alpha} - \eps}(\R^N) \cap L^{\frac{N}{\alpha}\frac{N+\alpha}{N-2s} + \eps}(\R^N)$ for some $\eps>0$; in particular, the convolution is well defined if $u \in L^2(\R^N) \cap L^{\frac{N}{\alpha} \frac{2N}{N-2s}}(\R^N)$, and thus if $u\in L^1(\R^N) \cap L^{\infty}(\R^N)$.
\end{Remark}

We recall now the definitions of weak solution, and of viscosity solution (see for instance \cite[page 136]{SeV} or \cite[Definition 2.1]{CFQ}).

\begin{Definition}[Weak solution]
Let $\Omega \subseteq \R^N$ and $g: \Omega 
\to \R$ be measurable. 
We say that $u\in H^s(\Omega)$ is a \emph{weak} subsolution [supersolution] of 
$$(-\Delta)^s u = g(x
) \quad \hbox{ in $\Omega$}$$
if
\begin{align}\label{eq_weak_sol}
&\int_{\R^N} (-\Delta)^{s/2} u (-\Delta)^{s/2} \varphi dx + \mu \int_{\R^N} u \varphi dx \leq \int_{\R^N} g(x
) \varphi dx
\\ &\qquad \qquad
\Big[\int_{\R^N} (-\Delta)^{s/2} u (-\Delta)^{s/2} \varphi dx + \mu \int_{\R^N} u \varphi dx \geq \, \int_{\R^N} g(x
) \varphi dx \Big] \notag
\end{align}
is well defined and holds for each nonnegative $\varphi \in X^s_0(\Omega)$. 
We say that $u$ is a weak solution if it is both a weak subsolution and a weak supersolution, i.e. if it satisfies the equality in \eqref{eq_weak_sol} for every $\varphi \in X^s_0(\Omega)$. 
Notice that, when $\Omega=\R^N$, we have $X^s_0(\R^N) \equiv H^s(\R^N)$.
\end{Definition}

\begin{Definition}[Viscosity solution] 
Let $\Omega \subseteq \R^N$ and $g: \Omega 
\to \R$. 
We say that $u\in C(\R^N)$ is a \emph{viscosity} subsolution [supersolution] of
$$(-\Delta)^s u = g(x
) \quad \hbox{ in $\Omega$}$$
 if, for any $x_0 \in \Omega$, every $U\subset \Omega$ open neighborhood of $x_0$, and every $\phi \in C^2(U)$ such that
$$\phi(x_0)=u(x_0), \quad \phi \geq u \; [\phi \leq u] \; \hbox{ in $U$}$$
set
$v:= \phi \chi_U + u \chi_{U^c}$ 
we have
\begin{align}\label{eq_visc_sol}
&(-\Delta)^s v(x_0) \leq g(x_0
) \\ &\qquad 
 \big[(-\Delta)^s v(x_0) \geq g(x_0
)\big]. \notag
\end{align}
We say that $u$ is a viscosity solution if it is both a viscosity subsolution and a viscosity supersolution.
\end{Definition}

We observe that, generally, the function $v$ appearing in the definition of viscosity solution might be discontinuous. 
More generally, this definition involves lower and upper semicontinuity of $u$ (see for instance \cite[Definition 2.2]{CafSil2}). 
Furthermore, one can easily check that every (continuous) classical solution is a viscosity solution, that the sum of two subsolutions is still a subsolution (with source the sum of the sources), and that the notion of subsolution is conserved on subdomains $\Omega' \subset \Omega$.

We refer to \cite[Remark 2.11]{ROS1} and \cite[Theorem 1]{SeV} for some discussions on the relation between classical, weak and viscosity solutions on bounded domains.

The above definitions apply, mutatis mutandis, to equation depending on $u$, i.e. where the right hand side is of the form $h_u(x)$; in this case (fixed $u$) the definition applies to $g(x):=h_u(x)$. In particular this adapts to our nonlocal equation
\begin{equation}\label{eq_main_prel}
(-\Delta)^s u + \mu u = (I_{\alpha}*F(u))f(u) \quad \hbox{on $\R^N$}
\end{equation}
by substituting \eqref{eq_weak_sol} with 
$$\int_{\R^N} (-\Delta)^{s/2} u (-\Delta)^{s/2} \varphi dx + \mu \int_{\R^N} u \varphi dx \leq \int_{\R^N} \big(I_{\alpha} * F(u)\big) f(u) \varphi dx,$$
where we implicitly assume \hyperref[(f1)]{\textnormal{(f1)}}-\hyperref[(f2)]{\textnormal{(f2)}} to give sense to the integrals, and substitute \eqref{eq_visc_sol} with
$$(-\Delta)^s v(x_0) \leq \big(I_{\alpha}*F(u)\big)(x_0) f(u(x_0)) ;$$
in this last case, we need some assumptions on $f$ and $u$ to have $I_{\alpha}*F(u)$ well defined pointwise, see Remark \ref{rem_conv_welldef}.

\medskip

In Appendix \ref{app_ex_comp} we collect some standard lemmas on existence results and comparison principles, both for weak and viscosity solutions.

\subsection{Fractional auxiliary functions}
\label{sec_frac_aux}

In order to implement some comparison argument, we search for a function which behaves like $\sim\frac{1}{|x|^{\beta}}$, $\beta >0$, and which lies in $H^s(\R^N)$: in order to handle the presence of a pole in the origin when $\beta \geq N$, we 
make the following choice, by considering, for any $\beta>0$, 
$$h_{\beta}(x):= \frac{1}{(1+|x|^2)^{\frac{\beta}{2}}};$$
notice that, when $\beta=N+2s$, this function is related to the extremals of the fractional Sobolev inequality \cite{Lie2} and to the solutions of the zero mass critical fractional Choquard equation \cite{Le0}. 
Chosen $h_{\beta}$ in this way, we have \cite[Table 1 page 168]{Kwa1} 
\begin{equation}\label{eq_frac_hyper}
(-\Delta)^s h_{\beta} (x) = C_{\beta, N, s} \, {}_2F_1\left(\frac{N}{2}+s, \frac{\beta}{2} + s, \frac{N}{2}; -|x|^2\right)
\end{equation}
where $C_{\beta, N, s}:= 2^{2s} \frac{\Gamma\big(\frac{N}{2}+s\big) \Gamma\big(\frac{\beta}{2} + s\big)}{\Gamma\big(\frac{N}{2}\big) \Gamma\big(\frac{\beta}{2}\big)}>0 $ 
and ${}_2F_1$ 
denotes the Gauss hypergeometric function (see also \cite[Corollary 2]{DKK}, observed that $h_{\beta}(x)= {}_2F_1(\frac{N}{2}, \frac{\beta}{2}, \frac{N}{2}, -|x|^2)$). 
Notice that we will be interested in 
$$\beta \in (0, N+2s].$$

In Appendix \ref{app_frac_dec} we collect some results on Gauss hypergeometric functions and their asymptotic behaviour at infinity. We use now this auxiliary function to study some comparison function.

\begin{Lemma}[Comparison for weak equation]\label{lem_confr_2xbeta}
Let $u\in C(\R^N)$ be a weak solution 
of
\begin{equation}\label{eq_comparison_hb}
(-\Delta)^s u + \lambda u = \gamma h_{\beta} \quad \hbox{in $\R^N\setminus B_{\rho}(0)$}
\end{equation}
for some $\lambda, \gamma>0$, $\rho>0$ and 
$\beta \in \left(\frac{N}{2}, N+2s\right].$
Then 
$$ \limsup_{|x|\to +\infty} u(x) |x|^{\beta} <\infty.$$
Moreover, if $\beta \in (\frac{N}{2}, N+2s)$, we have
$$\lim_{|x|\to +\infty} u(x) |x|^{\beta} = \frac{\gamma}{\lambda}.$$
\end{Lemma}

\claim Proof.
We start noticing that, since $\beta > \frac{N}{2}$, then the equation is well posed from a weak point of view.
By \cite[Lemma A.3]{CG0} (see also \cite[Lemma 4.3]{FQT}) there exists a continuous function $w \in H^{2s}(\R^N)$, such that
$$(-\Delta)^s w + \lambda w =0 \quad \hbox{in $\R^N\setminus B_{\rho}(0)$}$$
in the weak sense and pointwise, and moreover, for some $C_1'',C_2''>0$,
$$\frac{C_1''}{|x|^{N+2s}} <w(x) \leq \frac{C_2''}{|x|^{N+2s}}, \quad \hbox{for every $|x|>\rho$}.$$
Let thus define, for some $\tau, \sigma \in \R$ and 
$\theta \in [\beta, N+2s]$ to be chosen,
$$v_{\tau,\sigma}(x):= \frac{\gamma}{\lambda}h_{\beta}(x)+ \sigma h_{\theta}(x) + \tau w(x)$$
for every $x \in \R^N$. 
We have, for $|x|>\rho$,
\begin{align*}
(-\Delta)^s v_{\tau,\sigma}(x) + \lambda v_{\tau,\sigma}(x) &= \gamma h_{\beta}(x) + \left( \frac{\gamma}{\lambda} (-\Delta)^s h_{\beta}(x) + \sigma (-\Delta)^s h_{\theta}(x) + \lambda \sigma h_{\theta}(x)\right) \\
&=: \gamma h_{\beta}(x) + g_{\sigma, \theta}(x).
\end{align*}
By Lemma \ref{lem_calcol_2potenz} we obtain
\begin{itemize}
\item if $\beta \in (\frac{N}{2}, N)\setminus \{ N-2s\}$,
$$g_{\sigma, \theta}(x) \sim \frac{\gamma}{\lambda} C'_{\beta, N, s} h_{\beta +2s}(x) + \sigma C'_{\theta, N, s} h_{\theta+2s}(x) + \lambda \sigma h_{\theta}(x) \quad \hbox{as $|x|\to +\infty$};$$
in this case we assume $\theta \in (\beta, \min\{N, \beta+2s\})\setminus \{N-2s\}$;

\item if $\beta =N$,
$$g_{\sigma, \theta}(x) \sim \frac{\gamma}{\lambda} C'_{N, N, s} \log(x) h_{N+2s}(x) + \sigma C'_{\theta, N, s} h_{N+2s}(x) + \lambda \sigma h_{\theta}(x) \quad \hbox{as $|x|\to +\infty$};$$
in this case we assume $\theta \in (N, N+2s)$;

\item otherwise 
$$g_{\sigma, \theta}(x) \sim \frac{\gamma}{\lambda} C'_{\beta, N, s} h_{N+2s}(x) + \sigma C'_{\theta, N, s} h_{N+2s}(x) + \lambda \sigma h_{\theta}(x) \quad \hbox{as $|x|\to +\infty$},$$
and in this case 
\begin{itemize}
\item if $\beta=N-2s$ (possible only if $N>4s$), we choose $\theta \in (N, N+2s)$;
\item if $\beta \in (N,N+2s)$, we choose $\theta \in (\beta, N+2s)$;
\item if $\beta =N+2s$, we simply assume $\theta = 
N+2s$.
\end{itemize}
\end{itemize}
Assume first $\beta < N+2s$. 
By the abovementioned choices of $\theta >\beta$ we obtain
$$g_{\sigma, \theta}(x) \sim \lambda \sigma h_{\theta}(x) \quad \hbox{as $|x|\to +\infty$}. $$
In particular, fixed $\eps>0$, for some $R =R_{\eps} (\gamma, \lambda, \beta, \theta, \sigma) \gg 0$ (we may assume $R>\rho$) we obtain
$$(1-\eps) \lambda \sigma h_{\theta}(x) \leq g_{\sigma, \theta}(x) \leq (1+\eps) \lambda \sigma h_{\theta}(x) \quad \hbox{for $|x|\geq R$}$$
if $\sigma >0$, and
$$(1+\eps) \lambda \sigma h_{\theta}(x) \leq g_{\sigma, \theta}(x) \leq (1-\eps) \lambda \sigma h_{\theta}(x) \quad \hbox{for $|x|\geq R$}$$
if $\sigma <0 $. Notice that $R$ does not depend on $\tau$. Thus
$$(-\Delta)^s v_{\tau,\overline{\sigma}}(x) + \lambda v_{\tau,\overline{\sigma}}(x) \geq \gamma h_{\beta}(x) + (1-\eps) \lambda \overline{\sigma} h_{\theta}(x) \geq \gamma h_{\beta}(x) \quad \hbox{in $\R^N\setminus B_R(0)$}$$
by choosing a whatever $\overline{\sigma}>0$, and
$$(-\Delta)^s v_{\tau,\underline{\sigma}}(x) + \lambda v_{\tau,\underline{\sigma}}(x) \leq \gamma h_{\beta}(x) + (1-\eps) \lambda \underline{\sigma} h_{\theta}(x) \leq \gamma h_{\beta}(x) \quad \hbox{in $\R^N\setminus B_R(0)$}$$
by choosing a whatever $\underline{\sigma}<0$. 
Summing up
\begin{equation}\label{eq_confr_sistem}
\parag{ 
& (-\Delta)^s v_{\tau,\overline{\sigma}}(x) + \lambda v_{\tau,\overline{\sigma}}(x) \geq \gamma h_{\beta}(x) & \quad \hbox{in $\R^N\setminus B_R(0)$}, \\ 
& (-\Delta)^s v_{\tau,\underline{\sigma}}(x) + \lambda v_{\tau,\underline{\sigma}}(x) \leq \gamma h_{\beta}(x) & \quad \hbox{in $\R^N\setminus B_R(0)$}. }
\end{equation}
We choose now $\overline{\tau}>0$ such that
\begin{equation*}\label{eq_confr_tau_sop}
v_{\overline{\tau},\overline{\sigma}} - u \geq 0 \quad \hbox{on $B_R(0)$}.
\end{equation*}
Indeed, we impose
$$\frac{\gamma}{\lambda}h_{\beta}(x)+ \overline{\sigma} h_{\theta}(x) + \tau w(x) \geq u(x) \quad \hbox{on $B_R(0)$}$$
that is
$$\tau w(x) \geq u(x) - \frac{\gamma}{\lambda}h_{\beta}(x) - \overline{\sigma} h_{\theta}(x) \quad \hbox{on $B_R(0)$}$$
which is satisfied if we impose (recall that $\overline{\sigma}>0$)
$$ \tau \min_{B_R} w \geq \max_{B_R} u - \frac{\gamma}{\lambda} h_{\beta}(R) \geq u(x) - \frac{\gamma}{\lambda}h_{\beta}(x)- \overline{\sigma} h_{\theta}(x) \quad \hbox{on $B_R(0)$}$$
that is 
$ \overline{\tau} \geq \frac{ \max_{B_R} u - \frac{\gamma}{\lambda} h_{\beta}(R)}{\min_{B_R} w }.$
Similarly, we choose $\underline{\tau} \in \R$ 
such that
\begin{equation*}\label{eq_confr_tau_sot}
v_{\underline{\tau},\underline{\sigma}} - u \leq 0 \quad \hbox{on $B_R(0)$},
\end{equation*}
given by $ \underline{\tau} \leq \frac{ \min_{B_R} u - \frac{\gamma}{\lambda} h_{\beta}(R) }{\max_{B_R} w }$. 
We notice that both the minimum and the maximum of $w$ in the ball are finite and strictly positive, since $w>0$ is continuous. Thus, summing up
\begin{equation}\label{eq_confr_tau}
\parag{ & v_{\overline{\tau},\overline{\sigma}} - u \geq 0& \quad \hbox{on $B_R(0)$}, \\ &v_{\underline{\tau},\underline{\sigma}} - u \leq 0& \quad \hbox{on $B_R(0)$}.}
\end{equation}
By joining \eqref{eq_confr_sistem} with the assumption on $u$, we obtain
\begin{equation}\label{eq_confr_u_sistem}
\parag{ & (-\Delta)^s (v_{\overline{\tau},\overline{\sigma}} -u)(x)+ \lambda (v_{\overline{\tau},\overline{\sigma}}-u)(x) \geq 0 & \quad \hbox{in $\R^N\setminus B_R(0)$}, \\ 
& (-\Delta)^s (v_{\underline{\tau},\underline{\sigma}}-u)(x) + \lambda (v_{\underline{\tau},\underline{\sigma}}-u)(x) \leq 0 & \quad \hbox{in $\R^N\setminus B_R(0)$}. }
\end{equation}
By the weak version of the Comparison Principle (Lemma \ref{lem_comp_prin}) we obtain
$$
\parag{ & v_{\overline{\tau},\overline{\sigma}} - u \geq 0& \quad \hbox{on $\R^N$}, \\ &v_{\underline{\tau},\underline{\sigma}} - u \leq 0& \quad \hbox{on $\R^N$}.}
$$
that is
$$ \frac{\gamma}{\lambda}h_{\beta}(x)+ \underline{\sigma} h_{\theta}(x) + \underline{\tau} w(x) \leq u(x) \leq \frac{\gamma}{\lambda}h_{\beta}(x)+ \overline{\sigma} h_{\theta}(x) + \overline{\tau} w(x) $$
and hence, by the assumption on $w$,
$$ \frac{\gamma}{\lambda}h_{\beta}(x)+ \underline{\sigma} h_{\theta}(x) + \underline{\tau} \frac{C_1''}{|x|^{N+2s}} \leq u(x) \leq \frac{\gamma}{\lambda}h_{\beta}(x)+ \overline{\sigma} h_{\theta}(x) + \overline{\tau}\frac{C_2''}{|x|^{N+2s}}$$
for each $x \in \R^N$, $x\neq 0$. Thus
$$ \frac{\gamma}{\lambda}\frac{|x|^{\beta}}{(1+|x|^2)^{\frac{\beta}{2}}}+ \underline{\sigma} \frac{|x|^{\beta}}{(1+|x|^2)^{\frac{\theta}{2}}} + \underline{\tau} \frac{C_1''}{|x|^{N+2s-\beta}} 
\leq u(x)|x|^{\beta} 
\leq \frac{\gamma}{\lambda} \frac{|x|^{\beta}}{(1+|x|^2)^{\frac{\beta}{2}}}+ \overline{\sigma} \frac{|x|^{\beta}}{(1+|x|^2)^{\frac{\theta}{2}}} + \overline{\tau}\frac{C_2''}{|x|^{N+2s-\beta}},$$
which gives the claim passing to the limit $|x|\to +\infty$, since $\theta>\beta$ and $N+2s>\beta$. 

\smallskip

Assume now $\beta=N+2s$, and choose $\theta = 
\beta=N+2s$. Now we have
$$g_{\sigma, \theta}(x) \sim \overline{C}_{\sigma} h_{N+2s}(x) \quad \hbox{as $|x| \to +\infty$}$$
where $\overline{C}_{\sigma}:=\frac{\gamma}{\lambda} C'_{N+2s,N,s} + \sigma C'_{N+2s, N, s}+ \lambda \sigma;$ 
recall that $C'_{N+2s,N,s} 
<0$.
We can choose proper $\overline{\sigma} \in \R$ such that $\overline{C}_{\overline{\sigma}}<0$, 
and thus the first equation in \eqref{eq_confr_sistem} still hold.
Since the sign of $\overline{\sigma}$ may be now different, 
we choose 
$ \overline{\tau} \geq \frac{ \max_{B_R} u - \frac{\gamma}{\lambda} h_{\beta}(R)- \min\{\overline{\sigma},0\}
}{\min_{B_R} w }.$ 
We come up then with the same proof, obtaining
$$ \limsup_{|x|\to +\infty} u(x)|x|^{\beta} \leq \frac{\gamma}{\lambda}+ \overline{\sigma} + \overline{\tau}C_2''.$$
Notice that the appearing constants depend on $u, \gamma,\lambda, \rho, \beta, N, s$.
\QED

\medskip

\begin{Lemma}[Comparison for pointwise equation]\label{lem_confr_3xbeta}
Let $u \in C(\R^N)$ be a pointwise solution 
of \eqref{eq_comparison_hb}.
Then the conclusions of Lemma \ref{lem_confr_2xbeta} holds.
%
\end{Lemma}

\claim Proof.
The proof goes as the previous Lemma, with the difference that at the end we apply the pointwise version of the Comparison Principle (Lemma \ref{lem_comp_3prin}).
\QED

\section{Some preliminary estimates}
\label{sec_prel_est}

We start with some observations.

\begin{Remark}\label{rem_dec_N}
Let $u\in L^q(\R^N)$, for some $q \in [1,+\infty)$, be continuous and such that $|u|$ is radially symmetric and decreasing. Then, for every $x \in \R^N$,
\begin{align*}
|u(|x|)|^q |x|^N &= N |u(|x|)| \int_0^{|x|} t^{N-1} dt = N \int_0^{|x|} |u(|x|)|^q t^{N-1} dt \\
&\leq N \int_0^{|x|} |u(t)|^q t^{N-1} dt = \frac{N}{\omega_{N-1}} \int_{B_{|x|}(0)} |u(y)|^q dy \leq \frac{N}{\omega_{N-1}} \norm{u}_{L^q(\R^N)}^q
\end{align*}
where $\omega_{N-1}$ denotes the area of the $N-1$ dimensional sphere. Thus 
$$|u(x)| \leq \frac{C^2_u}{|x|^{\frac{N}{q}}}, \quad x \neq 0$$
where $C^2_u := C_N \norm{u}_q^q >0$. In particular, if $u \in L^1(\R^N)$, we have
$$|u(x)| \leq \frac{C^2_u}{|x|^N}, \quad x \neq 0.$$
\end{Remark}

\medskip

We keep with some preliminary lemmas; see \cite[Lemma 6.2]{MS0} (and \cite[Lemma C.3]{FLS}) for the first.
%
%
\begin{Lemma}[\cite{MS0}]\label{lem_stima_1_Riesz}
Let $g \in L^{\infty}(\R^N)$ continuous and $\theta>N$ be such that
$$\sup_{x \in \R^N} |g(x)| |x|^{\theta} < +\infty.$$
Then there exists $C=C(N,\alpha)>0$ such that
$$\pabs{\int_{\R^N} \frac{g(y)}{|x-y|^{N-\alpha}} dy - \frac{1}{|x|^{N-\alpha}} \int_{\R^N} g(y) dy} \leq \frac{C \norm{g}_{\infty,\theta}}{|x|^{N-\alpha}} \left( \frac{1}{1+|x|} + \frac{1}{1+ |x|^{\theta-N}}\right)$$
for each $x \in \R^N$, $x \neq 0$, where we recall that $\norm{g}_{\infty, \theta}= \norm{g(\cdot )(1+|\cdot|^{\theta})}_{\infty}$.
\end{Lemma}



\begin{Lemma}\label{lem_decay_Ialpha}
Let $u\in L^1 (\R^N)$ continuous be such that $|u|$ is radially symmetric and decreasing. Let $f$ satisfy \hyperref[(f1)]{\textnormal{(f1)}} and \hyperref[(f2)]{\textnormal{(f2,i)}}, and let $\theta \in (N, N+\alpha ]$. Then there exists $C=C(N,\alpha)>0$ such that
$$\pabs{\big(I_{\alpha}*F(u)\big)(x) - I_{\alpha}(x) \int_{\R^N} F(u)} \leq C \norm{F(u)}_{\infty,\theta} I_{\alpha}(x) \left( \frac{1}{1+|x|} + \frac{1}{1+ |x|^{\theta-N}}\right)$$
for each $x \in \R^N$, $x \neq 0$.
\end{Lemma}

\claim Proof.
First notice that $u\in L^{\infty}(\R^N)$, $F(u) \in L^{\infty}(\R^N)$, and that $I_{\alpha}*F(u)$ and $\int_{\R^N} F(u) $ are finite and well defined. 
By Remark \ref{rem_dec_N} we have
$$|u(x)| \leq \frac{C_u^2}{|x|^N} \to 0.$$
Thus $\big|F(u(x))\big| |x|^{\theta}$ is bounded on a ball $B_R$ (since $F(u)$ is bounded), and it is bounded on the complement of this ball since
$$\big|F(u(x))\big| |x|^{\theta} = \frac{\big|F(u(x))\big|}{|u(x)|^{\frac{N+\alpha}{N}}} |u(x)|^{\frac{N+\alpha}{N}} |x|^{\theta} \leq \frac{\big|F(u(x))\big|}{|u(x)|^{\frac{N+\alpha}{N}}} \frac{C}{|x|^{N+\alpha -\theta}}$$
by considering the growth condition \hyperref[(f2)]{\textnormal{(f2,i)}} of $F$ in zero (when $R\gg 0$, not depending on $\theta$) and the restriction on $\theta$. Thus
$\sup_{x \in \R^N} \big|F(u(x))\big| |x|^{\theta} < +\infty$ 
and Lemma \ref{lem_stima_1_Riesz} applies with $g(x):= F(u(x))$, which concludes the proof. We further notice 
 that
$$\norm{F(u)}_{\infty,\theta}\leq \norm{F(u)}_{\infty}(1+ R^{\theta}) + \left(\limsup_{t \to 0} \frac{|F(t)|}{|t|^{\frac{N+\alpha}{N}}} \right)\frac{1+R^{\theta}}{R^{N+\alpha}}$$
for any $\theta \in (N, N+\alpha]$ and any $R\gg 0$ (not depending on $\theta$, but depending on $u$).
\QED

\begin{Remark}\label{rem_cond_asym_dec}
 In what follows, for the sake of exposition we will restrict our analysis to the space of radially symmetric and decreasing functions in $L^1(\R^N)$, but we highlight that this assumption is needed only to get the a priori asymptotic decay of Remark \ref{rem_dec_N}. 
By the above proof, actually we see that we may ask only
$$|u(x)| \leq \frac{C}{|x|^{\omega}}$$
for some $\omega$ such that
$$\omega > \frac{N^2}{N+\alpha}.$$
In particular $\omega=N$, obtained in Remark \ref{rem_dec_N}, fits this condition. 
Alternatively, one may assume this a priori asymptotic decay on $u$ (and adapt the restrictions on $\theta$ by $\theta \in (N, \frac{N+\alpha}{N} \omega]$).
\end{Remark}

\begin{Corollary}\label{corol_stima_Riesz}
Let $u\in L^1 (\R^N)$ continuous be such that $|u|$ is radially symmetric and decreasing. 
Let $f$ satisfy \hyperref[(f1)]{\textnormal{(f1)}} and \hyperref[(f2)]{\textnormal{(f2,i)}}, and let $\theta \in (N, N+\alpha]$. 
Then for any $\eps>0$, there exists $R_{\eps}=R_{\eps}(N, \alpha,\theta)\gg 0 $ such that
$$\Big|\big(I_{\alpha}*F(u)\big)(x) \Big| \leq I_{\alpha}(x) \left( \pabs{\int_{\R^N} F(u)} + \eps \norm{F(u)}_{\infty,\theta}\right)$$
and
$$\big(I_{\alpha}*F(u)\big)(x) \geq I_{\alpha}(x) \left(\int_{\R^N} F(u)-\eps \norm{F(u)}_{\infty, \theta}\right) $$
for each $|x|\geq R_{\eps}$.
\end{Corollary}

\begin{Remark}\label{rem_dec_2sN}
In \cite{CGT3} it was showed that the solutions decay as fast as $\sim \frac{1}{|x|^{N+2s}}$ when the nonlinearity is linear or superlinear. 
In the sublinear case, we expect a slower decay. 
Indeed, assume \hyperref[(f1)]{\textnormal{(f1)}}, \hyperref[(f2)]{\textnormal{(f2)}} and \hyperref[(f4)]{\textnormal{(f4)}}, and let $u$ be a strictly  positive 
solution of \eqref{eq_main_prel}.
By \cite[Lemma 5.3]{CGT3} we have $u(x)\to 0$ as $|x|\to +\infty$, %
thus there exists $R\gg 0$ such that $0 \leq u(x) \leq \delta <1$ for $|x|\geq R$ and hence, by \eqref{eq_cond_4sublin_f}, 
$$f(u(x))\geq \underline{C} u^{r-1}(x) \quad \hbox{for $|x|\geq R$}$$
together with
$$u^{r-1}(x) \geq u(x) \quad \hbox{ for $|x|\geq R$}.$$
If we assume $(I_{\alpha}*F(u))(x) \geq 0$ for $|x|\geq R$, we gain 
$$(-\Delta)^s u + \mu u \geq (I_{\alpha}*F(u))u \quad \hbox{on $\R^N\setminus B_R(0)$}$$
which implies
$$(-\Delta)^s u + \tfrac{3}{2}\mu u \geq \left(I_{\alpha}*F(u) + \tfrac{1}{2}\mu\right) u \quad \hbox{on $\R^N\setminus B_R(0)$}.$$
By \cite[Proposition 4.5]{CGT3} (see also see also Remark \ref{rem_conv_welldef}) we have that $(I_{\alpha}*F(u))(x)\to 0$ as $|x|\to +\infty$, %
thus for some $R'\geq R \gg0$ we have
$$(-\Delta)^s u + \tfrac{3}{2}\mu u \geq 0 \quad \hbox{on $\R^N\setminus B_{R'}(0)$}.$$
At this point (being $u$ strictly positive) we conclude as in the proof of \cite[Theorem 1.3]{CGT3} and obtain
$$u(x) \geq \frac{C^1_u}{|x|^{N+2s}} \quad \hbox{for $|x| \geq R$}$$ 
for some constant $C^1_u = C_{N,\alpha, R, \mu} \min_{B_R} u>0$ and some sufficiently large $R\gg0$. 
\end{Remark}

By Remarks \ref{rem_dec_2sN} and \ref{rem_dec_N}, we obtain that every strictly positive, 
continuous, radially symmetric and decreasing solution of \eqref{eq_main_prel} in $L^1(\R^N)$ satisfies
\begin{equation}\label{eq_stima_grezza}
\frac{C^1_u}{|x|^{N+2s}} \leq u(x) \leq \frac{C^2_u}{|x|^N} \quad \hbox{for $|x|\geq R \gg 0$},
\end{equation}
 whenever $f$ satisfies \hyperref[(f1)]{\textnormal{(f1)}}-\hyperref[(f2)]{\textnormal{(f2)}} and \hyperref[(f4)]{\textnormal{(f4)}}, together with $\int_{\R^N} F(u)>0$: indeed in this case, by Lemma \ref{lem_decay_Ialpha}, we have $\big(I_{\alpha}*F(u)\big)(x) \sim I_{\alpha}(x) \int_{\R^N} F(u) >0$ for $|x|$ large.
Thus the goal is to improve the asymptotic decay \eqref{eq_stima_grezza} in the case of sublinear nonlinearities.

\medskip 

We highlight that, by Lemma \ref{lem_confr_2xbeta}, Corollary \ref{corol_stima_Riesz}, and a bootstrap argument one can give a first qualitative (not rigorous) proof of the main result. We refer to \cite[Remark 4.6.22]{GalT} for details.

\section{Estimate from above}
\label{sec_above}

First, we deal with the estimate from above. In this case we succeed in arguing in the weak sense with no additional assumption on $f$. 
In what follows we notice that, when $r> \frac{N+\alpha}{N}$, we are actually improving \eqref{eq_stima_grezza}.

\begin{Proposition}\label{prop_estim_above}
Assume \hyperref[(f1)]{\textnormal{(f1)}} and \hyperref[(f3)]{\textnormal{(f3)}}. 
Let $u\in H^s(\R^N)\cap L^1(\R^N)$, continuous, nonnegative, 
radially symmetric and decreasing, be a weak solution of \eqref{eq_main_prel}. 
Assume moreover
$$\mu> (r-1) \bar{C}^{\frac{1}{r-1}}.$$
Then, set $\beta$ as in \eqref{eq_defin_beta}
we have, for some $C_u\geq 0$, 
\begin{equation}\label{eq_stima_alto_prop}
\limsup_{|x|\to +\infty} u(x) |x|^{\beta} \leq C_u;
\end{equation}
if $\beta< N+2s$, the constant $C_u$ depends on $u$ in the following way:
$$C_u:= \frac{(2-r) \left(C_{N,\alpha} \pabs{ \int_{\R^N} F(u)}\right)^{\frac{1}{2-r}} }{\mu - (r-1)\bar{C}^{\frac{1}{r-1}} }
$$
where $C_{N,\alpha} >0$ is given in \eqref{eq_def_Riesz}.
\end{Proposition}

%
%
%

\claim Proof.
We start noticing that, by the Young product inequality, we obtain 
$$(I_{\alpha}*F(u))f(u) \leq \frac{1}{a} \big|I_{\alpha}*F(u)\big|^a + \frac{1}{b} |f(u)|^b$$
when $a,b>0$, $\frac{1}{a} + \frac{1}{b}=1$. In particular we choose $b=\frac{1}{r-1}$ and thus $a= \frac{1}{2-r}>0$ (possible thanks to the sublinearity restriction on $r$); with this choice, by \eqref{eq_cond_2sublin_f} and the fact that $u(x)\to 0$ as $|x|\to +\infty$, we obtain
$$(I_{\alpha}*F(u))f(u) \leq (2-r) \big|I_{\alpha}*F(u)\big|^{\frac{1}{2-r}} + (r-1)\bar{C}^{\frac{1}{r-1}} u$$
for $|x|\geq R$, where $R=R(u)\gg 0$ is sufficiently large. 
By Corollary \ref{corol_stima_Riesz}, for a whatever fixed $\theta \in (N, N+\alpha]$ and any $\eps>0$
we obtain
\begin{align*}
(I_{\alpha}*F(u))f(u) &\leq (2-r) \left(I_{\alpha}(x) \left( \pabs{ \int_{\R^N} F(u)} + \eps \norm{F(u)}_{\infty, \theta}\right)\right)^{\frac{1}{2-r}} + (r-1)\bar{C}^{\frac{1}{r-1}} u \\
&= (2-r) C_{N,\alpha}^{\frac{1}{2-r}} \left( \pabs{ \int_{\R^N} F(u)}+ \eps \norm{F(u)}_{\infty, \theta}\right)^{\frac{1}{2-r}} \frac{1}{|x|^{\frac{N-\alpha}{2-r}}} + (r-1)\bar{C}^{\frac{1}{r-1}} u
\end{align*}
for every $|x| \geq R_{\eps}=R_{\eps}(u,N, \alpha, \theta)$, 
thus
$$(-\Delta)^s u + \mu u \leq (2-r) C_{N,\alpha}^{\frac{1}{2-r}} \left( \pabs{ \int_{\R^N} F(u)}+ \eps \norm{F(u)}_{\infty,\theta}\right)^{\frac{1}{2-r}} \frac{1}{|x|^{\frac{N-\alpha}{2-r}}} + (r-1)\bar{C}^{\frac{1}{r-1}} u .$$
Notice that $F(u) \nequiv 0$ (otherwise, by the equation, $u\equiv 0$ and the claim is trivial), thus we set
$$\gamma_{u, \eps}:= (2-r) C_{N,\alpha}^{\frac{1}{2-r}} \left( \pabs{ \int_{\R^N} F(u)} + \eps \norm{F(u)}_{\infty, \theta}\right)^{\frac{1}{2-r}} 
>
0$$
and 
$\lambda:=\mu - (r-1)\bar{C}^{\frac{1}{r-1}} >0$ 
we obtain
$$(-\Delta)^s u + \lambda u \leq \frac{\gamma_{u,\eps}}{|x|^{\beta}} \quad \hbox{ in $\R^N \setminus B_{R_{\eps}}(0)$};$$
notice that we use the fact that $\frac{1}{|x|^{\frac{N-\alpha}{2-r}}} \leq \frac{1}{|x|^{\beta}}$ for $|x|$ large.
For each $\delta>1$ we have $\frac{1}{|x|^{\beta}} \leq \delta h_{\beta}(x)$ when $|x|> R_{\delta} :=(\delta^{\frac{2}{\beta}}-1)^{-\frac{1}{2}}$; we may choose $R_{\delta,\eps}>\max\{R_{\delta}, R_{\eps}\}$. Thus
\begin{equation}\label{eq_u_2gammau}
(-\Delta)^s u + \lambda u \leq \delta \gamma_{u,\eps} h_{\beta}(x) \quad \hbox{ in $\R^N\setminus B_{R_{\delta, \eps}}(0)$}.
\end{equation}
We have $ h_{\beta} \in L^2(\R^N)$, since $2\frac{N-\alpha}{2-r}>N$.
By Lemma \ref{lem_lax_milgr}, being $u \in H^s(\R^N)$, 
there exists $v\in H^s(\R^N)$ such that
$$\parag{ &(- \Delta)^s v + \lambda v = \delta \gamma_{u,\eps} h_{\beta}(x)& \quad \hbox{ in $\R^N \setminus B_{R_{\delta,\eps}}(0)$}, \\ &v = u& \quad \hbox{on $B_{R_{\delta,\eps}}(0)$}.}$$
Joining the first equation with \eqref{eq_u_2gammau} we obtain
$$(-\Delta)^s (u-v) + \lambda(u-v) \leq 0 \quad \hbox{ in $\R^N \setminus B_{R_{\delta,\eps}}(0)$}$$
and thus, by the weak version of the Comparison Principle (Lemma \ref{lem_comp_prin}) we obtain
\begin{equation}\label{eq_stima_uv_2comp}
u \leq v \quad \hbox{on $\R^N$}.
\end{equation}
By Lemma \ref{lem_confr_2xbeta}, if $\beta<N+2s$, we can estimate $v$ by
$\limsup_{|x|\to +\infty} v(x)|x|^{\beta} \leq \frac{ \delta \gamma_{u,\eps}}{\lambda}.$ 
This relation, combined with \eqref{eq_stima_uv_2comp}, gives
$$\limsup_{|x|\to +\infty} u(x) |x|^{\beta} \leq \frac{ \delta \gamma_{u,\eps}}{\lambda}$$
for each $\delta>1$. In particular, as $\delta \to 1^+$ and $\eps \to 0^+$, we obtain the claim. 
If $\beta=N+2s$, we argue similarly (without moving $\delta$ and $\eps$).
\QED

\bigskip

Notice that, if we assume \eqref{eq_strongf3}, then one can choose every $\bar{C}>0$, and thus in particular every $\mu>0$ is allowed (see Remark \ref{rem_eps_c_zero}). Anyway, in the proof of Theorem \ref{thm_main}, we will see how to drop the restriction on $\mu$. 

%
%

\medskip

We observe that the previous estimate from above is still valid by considering viscosity solutions $u \in L^1(\R^N) \cap C(\R^N)$, see Section \ref{sec_est_bel_vis}. We leave the details to the reader.

\section{Fractional concave Chain Rule and estimate from below} \label{sec_est_bel_vis}

Next, we deal with the estimate from below. 
We need first some preliminary results, in order to deal with the fractional Laplacian of the concave power of a function: since it might happen that $u^{\theta} \notin H^s(\R^N)$ when $u \in H^s(\R^N)$ and $\theta \in (0,1)$, the weak formulation seems not to be appropriate. 
Similarly, $(-\Delta)^s u^{\theta}$ might be not well defined pointwise, even if $u$ is regular enough. 
Notice that knowing a priori that $u$ is continuous, radially symmetric and decreasing seems of no use. 
The idea is thus to treat the problem via viscosity formulation.

The following lemma is a well known result in the case of convex and Lipschitz functions (see \cite[Theorem 1.1]{CafSir}, \cite[Theorem 19.1]{Gar0}). 
We state it here in the case of concave (not globally Lipschitz) function, in the framework of viscosity solutions. Notice that we do not require $u$ to be in $L^2(\R^N)$.

\begin{Lemma}[C\'ordoba-C\'ordoba chain rule inequality]\label{lem_chain_rule}
Let $\varphi: I \to \R$ be a concave function, $I\subseteq \R$ interval, such that $\varphi \in C^1(I)$. 
Let $u: \R^N \to I$. 
\begin{itemize}
\item Let $\Omega \subset \R^N$, and assume $\varphi \in Lip(u(\Omega))$. Then
$$[\varphi(u)]_{H^s(\Omega)} \leq  \norm{\varphi'}_{L^{\infty}(u(\Omega))} [u]_{H^s(\Omega)}.$$
In particular, if $\varphi \in Lip(I)$ and $(-\Delta)^{s/2} u \in L^2(\R^N)$, then $(-\Delta)^{s/2} \varphi(u)\in L^2(\R^N)$ and
$$\norm{(-\Delta)^{s/2}\varphi(u)}_2 \leq \norm{\varphi'}_{L^{\infty}(I)} \norm{(-\Delta)^{s/2} u}_2. $$
\item If $u$ is defined pointwise, then
$$(-\Delta)^s (\varphi(u))(x) \geq \varphi'(u(x)) (-\Delta)^s u(x)$$
for every $x \in \R^N$ such that $(-\Delta)^s (\varphi(u))(x)$ and $(-\Delta)^s u(x)$ are well defined.
\item Assume in addition $\varphi$ invertible, increasing, with $\varphi^{-1} \in C^2$ increasing. 
If $u$ is a continuous viscosity supersolution of
$$(-\Delta)^s u \geq g \quad \hbox{ in $\Omega$}$$
for some function $g$ and $\Omega \subseteq \R^N$, then $\varphi(u)$ is a viscosity supersolution of
$$(-\Delta)^s (\varphi(u)) \geq \varphi'(u) g \quad \hbox{ in $\Omega$}.$$
\end{itemize}
\end{Lemma}

\claim Proof.
The first claim is a direct consequence of the Lipschitz continuity
$$ \int_{\Omega}\int_{\Omega} \frac{|\varphi(u(x))-\varphi(u(y))|^2}{|x-y|^{N+2s}} dx dy  \leq  \norm{\varphi'}_{L^{\infty}(u(\Omega))}^2  \int_{\Omega} \int_{\Omega} \frac{|u(x)-u(y)|^2}{|x-y|^{N+2s}} dx dy.$$
Secondly, by the concavity of $\varphi$, for each $t,r \in I$ we have
$$\varphi(t) - \varphi(r) \geq \varphi'(t) (t-r)$$
thus
\begin{align*}
(-\Delta)^s(\varphi(u))(x) &= C_{N,s} \int_{\R^N} \frac{\varphi(u(x))-\varphi(u(y))}{|x-y|^{N+2s}} dy \\
&\geq C_{N,s} \int_{\R^N} \frac{\varphi'(u(x)) \big(u(x)-u(y)\big)}{|x-y|^{N+2s}} dy = \varphi'(u(x)) (-\Delta)^s u(x).
\end{align*}

We move to the third part. Let $x_0 \in U \subset \Omega$ and $\phi \in C^2(U)$ be such that $\phi(x_0)= \varphi(u(x_0))$ and $\phi\leq \varphi(u)$ in $U$, and set $v:= \phi \chi_{U} + \varphi(u) \chi_{U^c}$. Let now
$$\psi:= \varphi^{-1} \circ \phi, \quad w:=\varphi^{-1} \circ v = \psi \chi_{U} + u \chi_{U^c}.$$
By the assumptions on $\varphi^{-1}$ we have $\psi \in C^2(U)$, $\psi(x_0)=u(x_0)$ and $\psi \leq u$ in $U$. Thus
$$(-\Delta)^s w(x_0) \geq g(x_0).$$
On the other hand, $w = \psi \in C^2$ on $U$ and $\varphi(w) = \phi \in C^2$ on $U$, hence both the functions are regular enough in a neighborhood of $x_0$ to state that both the fractional Laplacians are well defined (see Proposition \ref{prop_well_posed}). 
 Thus we may apply the previous point and obtain
$$(-\Delta)^s (\varphi(w))(x_0) \geq \varphi'(w(x_0)) (-\Delta)^s w(x_0).$$
Since $w(x_0)=u(x_0)$, $\varphi(w)=v$ and $\varphi'$ is positive, we obtain, by joining the two previous inequalities
$$(-\Delta)^s v (x_0) \geq \varphi'(u(x_0)) g(x_0)$$
which is the claim. 
This concludes the proof. 
\QED

\bigskip

As a corrollary, 
we obtain the following result.

\begin{Corollary}\label{corol_concav_u2}
Let $\theta \in (0,1)$, and let $u \in C(\R^N)$ be strictly positive. 
\begin{itemize}
\item 
We have
$$[u^{\theta}]_{H^s(\Omega)} \leq \frac{\theta}{\min_{\Omega} u^{1-\theta} } [u]_{H^s(\Omega)}
$$
for each $\Omega \subset \subset \R^N$. 
In particular, if $u \in H^s_{loc}(\R^N)$, 
then $ u^{\theta} \in H^s_{loc}(\R^N)$. %
As a consequence, if $u \in H^s(\R^N)$, then
$$[u^{\theta}]_{H^s(\Omega)} \leq \frac{\theta}{\min_{\Omega} u^{1-\theta} }\norm{(-\Delta)^{s/2} u}_2.
$$
\item If $(-\Delta)^s u$ is well defined pointwise, then
$$(-\Delta)^s u^{\theta}(x) \geq \frac{\theta}{(u(x))^{1-\theta}} (-\Delta)^s u(x)$$
for every $x \in \R^N$ such that $(-\Delta)^s u^{\theta}(x)$ is well defined. 
\item If $u$ is a viscosity supersolution of
$$(-\Delta)^s u \geq g \quad \hbox{ in $\Omega$}$$
for some function $g$ and $\Omega \subseteq \R^N$, then $u^{\theta}$ is a viscosity supersolution of
$$(-\Delta)^s u^{\theta} \geq \frac{\theta}{u^{1-\theta}} g \quad \hbox{ in $\Omega$}.$$
\end{itemize}
\end{Corollary}

We are now ready the prove the estimate from below.

\begin{Proposition}\label{prop_below_2}
Assume \hyperref[(f1)]{\textnormal{(f1)}}-\hyperref[(f2)]{\textnormal{(f2,i)}} and the sublinear condition \hyperref[(f4)]{\textnormal{(f4)}}. 
Let $u\in L^1(\R^N) \cap C(\R^N)$, strictly positive, radially symmetric and decreasing, be a viscosity solution of \eqref{eq_main_prel}. Assume $\int_{\R^N}F(u)>0$.
Then, 
$$\liminf_{|x|\to +\infty} u(x) |x|^{\frac{N-\alpha}{2-r}} \geq C'_u$$
where 
$$C'_u := \left( \frac{\underline{C} C_{N,\alpha} \int_{\R^N} F(u)}{\mu } \right)^{\frac{1}{2-r}} $$
and $C_{N,\alpha}>0$ is given in \eqref{eq_def_Riesz}. 
Moreover, set $\beta$ as in \eqref{eq_defin_beta},
we have, for some $C''_u>0$, 
$$\liminf_{|x|\to +\infty} u(x) |x|^{\beta} \geq C''_u;$$
if $\frac{N-\alpha}{2-r}\leq N+2s$ (i.e. $\beta=\frac{N-\alpha}{2-r}$), we have $C''_u:=C'_u$, otherwise we have $C''_u:=C^1_u$ (see Remark \ref{rem_dec_2sN}).

The result in particular applies to pointwise solutions. 
\end{Proposition}

\claim Proof.
First notice that, by the assumptions, $u\in L^1(\R^N) \cap L^{\infty}(\R^N)$ and thus, by Remark \ref{rem_conv_welldef}, $I_{\alpha}*F(u)$ is pointwise well defined.
By Corollary \ref{corol_concav_u2}, since $2-r \in (0,1-\frac{\alpha}{N}]\subset (0,1)$ we have 
$$(-\Delta)^s u^{2-r} \geq \frac{2-r}{u^{r-1}} \Big(-\mu u + \big(I_{\alpha}*F(u))f(u)\Big) $$
on $\R^N$, in the viscosity sense.
Thus
$$ (-\Delta)^s u^{2-r} + \mu (2-r) u^{2-r} \geq (2-r) \frac{\big(I_{\alpha}*F(u))f(u)}{u^{r-1}} .$$
For a fixed $\theta \in (N, N+\alpha]$ and any $\eps >0$ small, by Corollary \ref{corol_stima_Riesz} and \eqref{eq_cond_4sublin_f} (since $u(x) \to 0$ as $|x|\to +\infty$, being $u$ decreasing and in $L^1(\R^N)$)
we obtain -- we use here that $\int_{\R^N} F(u)>0$ -- 
$$\big(I_{\alpha}*F(u)\big)f(u) \geq \underline{C} \left(\int_{\R^N} F(u)-\eps \norm{F(u)}_{\infty, \theta}\right)
I_{\alpha} u^{r-1} \quad \hbox{in $\R^N \setminus B_{R_{\eps}}(0)$}$$
for some $R_{\eps} \gg 0$, thus
$$ (-\Delta)^s u^{2-r} + \mu (2-r) u^{2-r} \geq \left( (2-r) \underline{C}
\left(\int_{\R^N} F(u)-\eps \norm{F(u)}_{\infty,\theta}\right)
\right) I_{\alpha} \quad \hbox{in $\R^N \setminus B_{R_{\eps}}(0)$};$$
that is, exploiting $\frac{1}{|x|^{N-\alpha}} \geq \frac{1}{(1+|x|^2)^{\frac{N-\alpha}{2}}}$, we get
$$(-\Delta)^s u^{2-r} + \lambda' u^{2-r} \geq \gamma'_{u,\eps} h_{N-\alpha} \quad \hbox{ in $\R^N \setminus B_{R_{\eps}}(0)$}$$
in the viscosity sense, where
$$\gamma'_{u,\eps} := (2-r) \underline{C} C_{N,\alpha} 
\left(\int_{\R^N} F(u)-\eps \norm{F(u)}_{\infty,\theta}\right) >0$$
and $\lambda':=\mu (2-r).$
We observe that $u^{2-r} \in L^{\infty}(B_R(0)) \cap C(B_{R_{\eps}}(0))$, while $h_{N-\alpha} \in L^{\infty}(\R^N) \cap C^{\sigma}_{loc}(\R^N)$
 (for any $\sigma$), thus by Lemma \ref{lem_lax_3milgr}, there exists $v\in C^{\omega}_{loc}(\R^N)$, for some $\omega>2s$ such that
$$\parag{ 
&(- \Delta)^s v + \lambda' v = \gamma'_{u,\eps} h_{N-\alpha} 
& \quad \hbox{ in $\R^N \setminus B_{R_{\eps}}(0)$}, \\ &v = u^{2-r} 
& \quad \hbox{on $B_{R_{\eps}}(0)$},}$$
pointwise. 
Thus
$$(-\Delta)^s (u^{2-r}-v) + \lambda' (u^{2-r}-v) \geq 0 \quad \hbox{ in $\R^N \setminus B_{R_{\eps}}(0)$}$$
in the viscosity sense, with
$$u^{2-r} - v \geq 0 \quad \hbox{on $B_{R_{\eps}}(0)$}.$$
Observe that, by Lemma \ref{lem_confr_3xbeta}, we have $v(x) \to 0$ as $|x|\to +\infty$.
Since $(u^{r-2}-v)(x) \to 0$ as $|x|\to +\infty$, by the viscosity version of the Comparison Principle (Lemma \ref{lem_comp_3prin}) we obtain
$$u^{2-r} \geq v \quad \hbox{on $\R^N$}.$$
By Lemma \ref{lem_confr_3xbeta} we gain
$$\liminf_{|x|\to +\infty} v(x) |x|^{N-\alpha} \geq \frac{\gamma'_{u,\eps}}{\lambda'} .$$
Combining the previous inequalities and 
sending $\eps \to 0^+$, we have the first claim. We conclude by adapting Remark \ref{rem_dec_2sN} to the viscosity case (notice that $u\in L^1(\R^N) \cap L^{\infty}(\R^N)$).
\QED

\bigskip

%
%

By the results in \cite{CGT3, CG1}, we gain sufficient conditions in order to state that a weak solution is a pointwise solution. 

\begin{Corollary}\label{cor_est_bel_wea}
Assume \hyperref[(f1)]{\textnormal{(f1)}}-\hyperref[(f2)]{\textnormal{(f2,i)}} and the sublinear condition \hyperref[(f4)]{\textnormal{(f4)}}. 
Let $u\in H^s(\R^N) \cap L^1(\R^N) \cap C(\R^N)$, strictly positive, radially symmetric and decreasing, be a weak solution
of \eqref{eq_main_prel}. Assume moreover that $f \in C^{0,\sigma}_{loc}(\R)$ for some $\sigma \in (0,1]$ and $\int_{\R^N}F(u)>0$. Then $u$ is a classical solution and 
 the conclusions of Proposition \ref{prop_below_2} hold.
\end{Corollary}

Notice that, by the sublinearity in zero, $\sigma$ can lie only in $(0, r-1]$. 
We conjecture anyway that the conclusion of Corollary \ref{cor_est_bel_wea} holds in more general cases, by assuming $f$ merely continuous.

\section{Proofs of the main theorems} 
\label{sec_proof_main_2}

We can sum up some of the results of the previous sections in the following.

\begin{Corollary}
Assume \hyperref[(f1)]{\textnormal{(f1)}}-\hyperref[(f2)]{\textnormal{(f2,i)}} and the sublinear conditions \hyperref[(f3)]{\textnormal{(f3)}}-\hyperref[(f4)]{\textnormal{(f4)}}, in particular
$$0< \liminf_{t \to 0} \frac{f(t)}{|t|^{r-1}} \leq \limsup_{t \to 0} \frac{f(t)}{|t|^{r-1}} \leq \bar{C} < +\infty.$$ 
Let $u\in H^s(\R^N) \cap L^1(\R^N) \cap C(\R^N)$, 
strictly positive, radially symmetric and decreasing, be a weak solution of \eqref{eq_main_prel}. 
Finally assume $\mu> (r-1) \bar{C}^{\frac{1}{r-1}}$, 
$f\in C^{0,\sigma}(\R) \; \hbox{ for some $\sigma \in (0,r-1]$}$, and $\int_{\R^N}F(u)>0$. 
Then we have
$$0 < \liminf_{|x| \to +\infty} u(x)|x|^{\beta} \leq \limsup_{|x|\to +\infty} u(x) |x|^{\beta} < +\infty$$
where $\beta$ is defined in \eqref{eq_defin_beta}.
\end{Corollary}

%
%
%

We can now conclude the proof of the main theorem.

\medskip

\claim Proof of Theorem \ref{thm_main}.
First, we show how to remove the restriction on $\mu$ in Proposition \ref{prop_estim_above}. Indeed, for any $\kappa>0$ we can write $\big(I_{\alpha}*F(u)\big) f(u) \equiv \big(I_{\alpha}*F_{\kappa}(u)\big) f_{\kappa}(u)$, where $f_{\kappa}:=\frac{1}{\kappa} f$ and $F_{\kappa}:=\kappa F$. We can thus rewrite \hyperref[(f3)]{\textnormal{(f3)}} as
$$|f_{\kappa} (t)|\leq \frac{1}{\kappa} \overline{C} t^{r-1} \quad \hbox{for $t \in (0,\delta)$}.$$
Since in Proposition \ref{prop_estim_above} we did not use that $F$ is the primitive of $f$ (in particular, we did not apply \hyperref[(f3)]{\textnormal{(f3)}} to $F$), fixed a whatever $\mu>0$ we can choose $\kappa$ such that
$$\mu > (r-1) \left(\frac{\overline{C}}{\kappa}\right)^{\frac{1}{r-1}} >0,$$
that is a large $\kappa$ given by $\kappa > \left(\frac{r-1}{\mu}\right)^{r-1} \overline{C}$, and obtain 
$$\limsup_{|x|\to +\infty} u(x) |x|^{\beta} \leq C_{u,\kappa}$$
where, if $\beta < N+2s$, 
\begin{equation*}\label{eq_constant_k}
C_{u,\kappa}:= \frac{(2-r) \left( C_{N,\alpha} \kappa \pabs{ \int_{\R^N} F(u)}\right)^{\frac{1}{2-r}} }{\mu -(r-1)\left(\frac{\bar{C}}{\kappa}\right)^{\frac{1}{r-1}} }.
\end{equation*}
We notice, as we expect, that as $\mu \to 0$ then $\kappa \to +\infty$ and $C_{u,\kappa}\to +\infty$, while $C'_u$ defined in Proposition \ref{prop_below_2} is invariant under $\kappa$-transformations.

\smallskip

We show now the sharp decay; indeed, we search for a $\kappa$ such that $C_{u,\kappa}=C'_u$. By a straightforward analysis of $g(\kappa):=C_{u,\kappa}-C'_u$, $\kappa> \left(\frac{r-1}{\mu}\right)^{r-1} \overline{C}$, we find a (unique, explicit) zero $\kappa^*$ (which actually is a point of minimum) if only if $\overline{C}=\underline{C}$, i.e. if $f$ is exactly a power near the origin. 

\smallskip

Now, by the results in \cite{CGT3}, we have that every positive solution is bounded, and every bounded solution is in $H^{2s}(\R^N)\cap C(\R^N) \cap L^1(\R^N)$. By the previous results we conclude the proof.
\QED

\medskip

\begin{Remark}\label{rem_eps_c_zero}
We notice that, when $r \in (\frac{N+\alpha}{N},2)$, by assuming 
$\limsup_{t \to 0^+} \frac{f(t)}{t^{r-1}} \in (0, +\infty)$ 
we obtain that
$ \limsup_{t \to 0^+} \frac{f(t)}{t^{r-\eps-1}} =0 $ for any $\eps>0$ (such that $r-\eps \in [ \frac{N+\alpha}{N},2)$). Thus we may extend the estimate from above of Proposition \ref{prop_estim_above} to a whatever $\mu>0$ also by paying the cost of a slower decay at infinity, that is
$$ \limsup_{|x|\to +\infty} u(x) |x|^{\beta_{\eps}} \leq C_{u,\eps}< \infty$$
where
$ \beta_{\eps}:= \min\left\{ \frac{N-\alpha}{2-r+\eps}, N+2s\right\}.$
If $\beta< N+2s$, $C_{u,\eps}:= \frac{(2-r+\eps) \left( C_{N,\alpha} \pabs{ \int_{\R^N} F(u)}\right)^{\frac{1}{2-r}} }{\mu}.$ 
 Clearly, when condition \eqref{eq_strongf3} holds, then the above statement holds for $\eps=0$.
\end{Remark}

%


Before concluding the proof of Corollary \ref{corol_main}, we make clear what we mean by Pohozaev minimum.

\begin{Definition}\label{def_pohozaev_min2}
Consider the functional $\mc{I}\in C^1(H^s(\R^N), \R)$ defined by
$$\mc{I}(u):= \frac{1}{2} \int_{\R^N} |(-\Delta)^{s/2}u|^2 dx + \frac{1}{2} \int_{\R^N} u^2 dx - \frac{1}{2} \int_{\R^N} \big(I_{\alpha}*F(u)\big)F(u) dx;$$
it is easy to see that weak solutions of \eqref{eq_main_prel} are critical points of $\mc{I}$. We say that $u$ is a \emph{Pohozaev minimum} of $\mc{I}$ if
$$\mc{I}(u) = \inf \left\{ \mc{I}(v) \mid v \in H^s(\R^N)\setminus \{0\}, \; \mc{P}(v)=0 \right\}$$
where
$$\mc{P}(u):=\frac{N-2s}{2} \int_{\R^N} |(-\Delta)^{s/2}u|^2 dx + \frac{N}{2} \mu \int_{\R^N} u^2 dx - \frac{N+\alpha}{2} \int_{\R^N} \big(I_{\alpha}*F(u)\big)F(u) dx.$$
By \cite{CGT3} we know that every Pohozaev minimum is a critical point of $\mc{I}$, and thus every Pohozaev minimum $u$ satisfies
$$\mc{I}(u) = \inf \left\{ \mc{I}(v) \mid v \in H^s(\R^N)\setminus \{0\}, \; \mc{P}(v)=0, \; \mc{I}'(v)=0 \right\}.$$
We say instead that $u$ is a \emph{ground state} for $\mc{I}$ if
$$\mc{I}(u)= \inf \left \{ \mc{I}(v) \mid v \in H^s(\R^N)\setminus \{0\}, \; \mc{I}'(v)=0\right\}.$$
\end{Definition}

We point out that, by \cite{CGT3}, we know that there exists a (radially symmetric) Pohozaev minimum if we further assume that $F$ is nontrivial and not critical, i.e. the limits in \eqref{eq_condiz_F} are zero. Moreover, the notion of Pohozaev minimum and ground state coincide whenever all the critical points satisfy the Pohozaev identity, fact known under some restriction on $s, \alpha$ and $f$, see \cite{CGT5}.

\medskip

\claim Proof of Corollary \ref{corol_main}.
By the results in \cite{CG1}, we have that every Pohozev minimum has strict constant sign -- e.g., strictly positive -- (if $f$ is odd or even, and H\"older continuous), and it is radially symmetric and decreasing (if in addition $f$ has constant sign on $(0,+\infty)$). 
Thus we conclude by the previous results.
\QED

\bigskip

%
%

All the previous theorems particularly apply to homogeneous nonlinearities

\medskip

\claim Proof of Theorem \ref{thm_homog0_ws}, Corollary \ref{corol_homog0_gs}. 
Theorem \ref{thm_homog0_ws} is a direct consequence of the above result, where formally $f(t)=\sqrt{r}|t|^{r-2} t$. 
By \cite[Theorems 3.2 and 4.2]{DSS1} we have that every ground state satisfies all the assumptions of the previous results; thus we have the claim of Corollary \ref{corol_homog0_gs}. 
\QED


 \appendix

\section{Appendices} 

\subsection{Existence and comparison results}
\label{app_ex_comp}

We collect here some results which are already known in literature, even if the author was not able to find a precise reference.
We start with an existence result (see also \cite[Corollary 1.15]{ShSp}).

\begin{Lemma}[Existence for weak solutions]\label{lem_lax_milgr}
Let $\Omega\subset \R^N$ be of class $C^{0,1}$ with bounded boundary\footnote{The result is still valid in a whatever $\Omega^c$ \emph{extension domain} (see \cite{DPV}).}, $\lambda>0$, $\psi \in H^{s}(\Omega^c)$, 
and $g \in L^q(\Omega)$, for some $q \in [\frac{2N}{N+2s}, 2]$.
Then there exists a (unique) function $v \in H^s(\R^N)$ such that 
$$\parag{ &(- \Delta)^s v + \lambda v = g& \quad \hbox{in $\Omega$}, \\ &v = \psi& \quad \hbox{on $\Omega^c$,}}$$
in the weak sense, that is $v \in X^s_0(\Omega)+\psi$ and 
$$\int_{\R^N} (-\Delta)^{s/2} v (-\Delta)^{s/2} \varphi dx + \lambda \int_{\R^N} v \varphi dx = \int_{\R^N} g \varphi dx$$
for every $\varphi \in X^s_0(\Omega)$. 
If moreover $g \in L^{q}_{loc}(\R^N)$ for some $q>\frac{N}{2s}$, then $v\in L^{\infty}_{loc}(\R^N)$. If instead $g \in C^{0, \sigma}_{loc}(\R^N)$ for some $\sigma \in (0,1]$, then $v \in C^{2s+\sigma}_{loc}(\R^N)$.
\end{Lemma}


\claim Proof.
By \cite[Theorem 5.4]{DPV} we know that there exists $\tilde{\psi} \in H^s(\R^N)$ such that $\tilde{\psi}_{|\Omega^c}\equiv \psi$.
The problem is thus equivalent to
$$\parag{ &(- \Delta)^s v + \lambda v = g& \quad \hbox{in $\Omega$}, \\ &v = \tilde{\psi}& \quad \hbox{on $\Omega^c$}.}$$
Consider $u= v-\tilde{\psi}$ and rewrite the weak formulation as
$$\int_{\R^N} (-\Delta)^{s/2} u (-\Delta)^{s/2} \varphi dx + \lambda \int_{\R^N} u \varphi dx = \int_{\R^N} (g-\lambda \psi) \varphi dx - \int_{\R^N} (-\Delta)^{s/2} \tilde{\psi} (-\Delta)^{s/2} \varphi dx.$$
It is easy to see that the left-hand side is a bilinear, continuous coercive map on $X^s_0(\Omega)$, while
$$\varphi \in X^s_0(\Omega) \mapsto \int_{\R^N} (g-\lambda \psi) \varphi dx - \int_{\R^N} (-\Delta)^{s/2} \tilde{\psi} (-\Delta)^{s/2} \varphi dx$$
belongs to the dual space $(X^s_0)^*(\Omega)$. By Lax-Milgram theorem, we obtain a solution $u \in X^s_0(\Omega)$, which implies that $v:=u+\tilde{\psi}$ is the desired function.

Finally, the regularity results are a consequence of De Giorgi-Nash-Moser estimates \cite[Proposition 2.6]{JLX} and Schauder estimates \cite[Theorem 2.11]{JLX}. 
\QED

\bigskip

The following existence result can be found in \cite[Lemma 2.2 and Remark 4.1]{CFQ} for bounded domains, and in \cite[Theorem A.1]{SoV} for the homogeneous case $\psi\equiv 0$.

\begin{Lemma}[Existence for viscosity solutions]\label{lem_lax_3milgr}
Let $\Omega \subset \R^N$ be a $C^2$ domain, $\lambda>0$, $\psi \in L^{\infty}(\Omega^c) \cap C(\Omega^c)$, 
and $g \in L^{\infty}(\Omega) \cap C(\overline{\Omega})$.
Then there exists a function $v \in C(\R^N) \cap L^{\infty}(\R^N)$
such that 
$$\parag{ &(- \Delta)^s v + \lambda v = g& \quad \hbox{in $\Omega$}, \\ &v = \psi& \quad \hbox{on $\Omega^c$,}}$$
in the viscosity sense. If $g \in C^{\sigma}_{loc}(\Omega)$ for some $\sigma \in (0,1)$, then $v\in C^{\gamma}_{loc}(\Omega)$, for some $\gamma>2s$ is a pointwise solution.
If $\psi \equiv 0$, we further have $v \in C^s(\R^N) \cap C^{\gamma}_{loc}(\Omega)$, for some $\gamma>\max\{1,2s\}$ and $\frac{w}{(\dist(\cdot,\partial \Omega))^s} \in C^{0,\theta}(\overline{\Omega})$ for some $\theta \in (0,1)$.
\end{Lemma}

\claim Proof.
First notice that, by extension, we may assume $g \in L^{\infty}(\R^N) \cap C(\R^N) \cap C^{\sigma}_{loc}(\Omega)$.
Since $\Omega$ is a $C^2$ domain, $g\in C(\R^N)$ and $\psi\in C(\Omega^c) \cap L^{\infty}(\Omega^c)$, by \cite[Theorem 4]{BCI} with $b\equiv c \equiv 0$ we obtain the existence of a (unique) viscosity solution $v \in C(\R^N)$, satisfying the boundary condition pointwise (see also \cite[page 615]{CafSil2}). 
Since the cited theorem is a corollary of \cite[Theorem 1]{BCI}, with $F(x,u,p,X,l)\equiv F(x,u,l)=l + \lambda u - g(x)$, $l=\mc{I}[u]\equiv (-\Delta)^s u$, $d \mu_x(z) = \frac{dz}{|z|^{N+\alpha}}$, one can notice, looking at the proof, that the found solution is actually bounded (see also \cite[Corollary 4]{SeV}). 
Thus $v$ is a bounded viscosity solution. 

By \cite[Theorem 2.6]{QX0}, since $(-\Delta)^s v =-\lambda v + g \in L^{\infty}(\Omega)$ 
with $v\in C(\overline{\Omega})$, we have $v \in C^{\gamma_1}_{loc}(\R^N)$ for some $\gamma_1>0$. 
Since $\psi \in L^{\infty}(\Omega^c)$ and $g-\lambda v \in C^{\min\{\sigma,\gamma_1\}}_{loc}(\Omega)$, by \cite[Theorem 2.5]{QX0} we have that $v\in C^{\gamma}_{loc}(\Omega)$ for some $\gamma>2s$; thus $(-\Delta)^s v$ is pointwise defined (actually H\"older continuous). 
As observed in \cite[Remark 2.3]{QX0}, we conclude that $v$ is a pointwise solution.
\QED

\bigskip

We write down the following two lemmas again for the reader's convenience. See \cite[Lemma A.1]{CG0} for the first.

\begin{Lemma}[Maximum Principle (weak)]\label{lem_comp_prin}
Let $\Omega \subset \R^N$ (possibly unbounded), $\lambda>0$, and let $u\in H^s(\R^N)$ be a weak subsolution of
$$(-\Delta)^s u + \lambda u \leq 0 \quad \hbox{in $\Omega$}.$$
Assume moreover that $u(x)\leq 0$ on $\Omega^c$. 
Then
\begin{equation}\label{eq_dis_comp_prin}
u(x)\leq 0 \quad \textit{on $\R^N$}.
\end{equation}
\end{Lemma}


\begin{Lemma}[Maximum Principle (viscosity)]\label{lem_comp_3prin}
Let $\Omega \subset \R^N$ open (possibly unbounded), $\lambda>0$, and let $u$ be a viscosity, continuous subsolution of
$$(-\Delta)^s u + \lambda u \leq 0 \quad \hbox{in $\Omega$}$$
such that
$$\lim_{|x|\to +\infty} u(x) \leq 0.$$
Assume moreover that $u(x)\leq 0$ on $\Omega^c$.
Then
\begin{equation}\label{eq_dis_comp_2prin}
u(x)\leq 0 \quad \textit{on $\R^N$}.
\end{equation}
The result applies, in particular, to pointwise solutions.
\end{Lemma}

\claim Proof.
We first observe that $u\in L^{\infty}(\R^N)$ and set $M:=\sup_{x \in \R^N} u(x)$. By contradiction, assume $M>0$. Let $(x_n)_n$ be a maximizing sequence, i.e. $u(x_n) \to M$ as $n\to +\infty$; we can assume that $x_n \in \Omega$. 
We observe that $(x_n)_n$ is bounded (up to a subsequence) since, if not, we would have $|x_n| \to +\infty$ and thus $\lim_n u(x_n) \leq 0$, 
which is an absurd. 
Thus $x_n \to x_0 \in \overline{\Omega}$, and by continuity $u(x_0)=M>0$; since $u(x)\leq 0$ on $\overline{\Omega^c} \supset \partial \Omega$, we have $x_0 \in \Omega$. In particular, $x_0$ is a point of maximum for $u$. 

We can thus choose a whatever $ U\subset \Omega $ neighborhood of $x_0$ and set $\phi \equiv u(x_0)$ as contact function in the definition of viscosity solution: indeed $\phi \in C^2(U)$, $\phi(x_0)=u(x_0)$ and $\phi \geq u$ in $U$. 
Hence, set $v:= \phi \chi_{U} + u \chi_{U^c}$ we have
\begin{align*}
0& \geq (-\Delta)^s v(x_0) + \lambda v(x_0) = C_{N,s} \int_{\R^N} \frac{u(x_0) - v(y)}{|x_0-y|^{N+2s}} dy + \lambda u(x_0) \\
&= C_{N,s} \int_{U^c} \frac{M - u(y)}{|x_0-y|^{N+2s}} dy + \lambda M >0,
\end{align*}
which is a contradiction. This concludes the proof.
\QED

\subsection{Decay of fractional auxiliary functions}
\label{app_frac_dec}

Recall that, by \eqref{eq_frac_hyper}, 
\begin{equation}\label{eq_frac_hyper_2}
(-\Delta)^s h_{\beta} (x) = C_{\beta, N, s} \, {}_2F_1\left(\frac{N}{2}+s, \frac{\beta}{2} + s, \frac{N}{2}; -|x|^2\right)
\end{equation}
where $h_{\beta}(x)= \frac{1}{(1+|x|^2)^{\frac{\beta}{2}}}$, and $C_{\beta, N, s}= 2^{2s} \frac{\Gamma\big(\frac{N}{2}+s\big) \Gamma\big(\frac{\beta}{2} + s\big)}{\Gamma\big(\frac{N}{2}\big) \Gamma\big(\frac{\beta}{2}\big)}>0 $.
We consider now the asymptotic behaviour at infinity of the hypergeometric function $ {}_2F_1 $ 
(see \cite[pages 559-560]{AS0}, but also \cite[pages 78-79, 88]{AAR} and \cite[page 161]{WG0}).
Recall that $\Gamma(z)$ is well defined whenever $z \in \R \setminus (-\N)$ and $|\Gamma(z)|\to +\infty$ as $z$ approaches $-\N$ (so that the \emph{reciprocal Gamma function} is well defined on $-\N$ and equals zero); 
moreover we have the symmetry property ${}_2F_1(a,b,c;x)={}_2F_1(b,a,c;x)$ and the fact that ${}_2F_1(0,b,c;x) =1$ and ${}_2F_1(-1,b,c;x) = 1 - \frac{b}{c} z$.
\begin{Lemma}[\cite{AS0}]
Consider ${}_2F_1(a,b,c; x)$. For the sake of simplicity, assume a priori that $a,b,c>0$ and
$$a-c \in \R_+ \setminus \N,$$
$$a-b \in \Z \iff a-b \in \N,$$ 
$$b-c \in \N \iff b-c \in \{0,1\};$$
in particular $a-b$ and $b-c$ do not lie in $\Z$ at the same time. 
We have the following asymptotic estimates as $x\to -\infty$.
\begin{itemize}
\item If $a-b \notin \Z$ and $b-c \notin \N$, then
$${}_2F_1(a,b,c; x) \sim \frac{\Gamma(c) \Gamma(b-a)}{\Gamma(c-a)\Gamma(b)} \frac{1}{(-x)^{a}} + \frac{\Gamma(c) \Gamma(a-b)}{\Gamma(c-b)\Gamma(a)} \frac{1}{(-x)^{b}} ;$$
\item If $b=c$ (and $a-b \notin \Z$), then
$${}_2F_1(a,b,b;x) =\frac{1}{(1-x)^{a}}; $$ 
\item If $b=c+1$ (and $a-b \notin \Z$), then
$${}_2F_1(a,b,b-1;x) = - \frac{\Gamma(b-1)\Gamma(b-a)}{\Gamma(b-a-1)\Gamma(b)}\frac{x}{(1-x)^{a+1}} + \frac{1}{(1-x)^{a+1}} \sim \frac{\Gamma(b-1)\Gamma(b-a)}{\Gamma(b-a-1)\Gamma(b)} 
 \frac{1}{(-x)^{a}};$$
\item If $a=b$ (and $b-c \notin \N$), then
$${}_2F_1(a,a,c;x) \sim \frac{\Gamma(c)}{\Gamma(a)\Gamma(c-a)} \frac{\log(-x)}{(-x)^a} + \frac{C_1}{(-x)^a} \sim \frac{\Gamma(c)}{\Gamma(a)\Gamma(c-a)} \frac{\log(-x)}{(-x)^a}; $$
\item If $a-b \in \N^*$ (and $b-c \notin \N$), then
$${}_2F_1(a,b,c;x) \sim \frac{\Gamma(c) \Gamma(a-b)}{\Gamma(c-b)\Gamma(a)} \frac{1}{(-x)^{b}} + C_2 \frac{\log(-x)}{(-x)^a} + \frac{C_3}{(-x)^a} \sim \frac{\Gamma(c) \Gamma(a-b)}{\Gamma(c-b)\Gamma(a)} \frac{1}{(-x)^{b}}. $$
\end{itemize}
Here $C_i$, $i=1,2,3$, are some strictly positive constants.
\end{Lemma}

Notice that $a= \frac{N}{2} + s$, $b=\frac{\beta}{2} + s$, $c= \frac{N}{2}$ satisfy the assumptions of the previous Lemma, whenever $s \in (0,1)$ and $\beta \in (0, N+2s]$. 
Thus, exploiting the representation of $(-\Delta)^s h_{\beta}$ given in \eqref{eq_frac_hyper_2} and the results on Gauss hypergeometric functions, we come up with the following estimates.

\begin{Lemma}\label{lem_calcol_2potenz}
Let $\beta \in (0,N+2s]$. 
Then $(-\Delta)^s h_{\beta}(x)$ is well-defined for every $x \neq 0$. Moreover, we have the following asymptotic behaviours:
\begin{itemize}
\item if $\beta\in (N, N+2s]$, then
$$(-\Delta)^s h_{\beta}(x) \sim C'_{\beta, N, s} \frac{1}{|x|^{N+2s}} \quad \hbox{ as $|x|\to +\infty$}$$
where $C'_{\beta,N,s}:=2^{2s} \frac{\Gamma\big(\frac{N}{2}+s\big) \Gamma\big(\frac{\beta}{2}-\frac{N}{2}\big)}{\Gamma\big(\frac{\beta}{2}\big) \Gamma\big(-s\big)}<0.$ 
This in particular includes the case $\beta = N-2s+2$ (possible if $s > 1/2$), with $C_{N-2s+2,N,s}'= - 2^{2s+1}\frac{s }{N-2s}<0$. 
Notice moreover that $C'_{N+2s,N,s}=2^{2s} \frac{\Gamma(s)}{\Gamma(-s)} \to 0$ as $s \to 1^{-}$.

\item if $\beta=N$, then 
$$(-\Delta)^s h_N(x) \sim C'_{N,N,s} \frac{\log(|x|)}{|x|^{N+2s}} \quad \hbox{ as $|x|\to +\infty$}$$
where 
$C'_{N,N,s}:=
2^{2s+1} \frac{\Gamma\big(\frac{N}{2}+s\big)}{\Gamma\big(\frac{N}{2}\big) \Gamma\big(-
s\big)} <0.$

\item if $\beta\in (N-2s, N)$, then
$$(-\Delta)^s h_{\beta}(x) \sim C_{\beta,N,s}' \frac{1}{|x|^{\beta + 2s}} \quad \hbox{as $|x|\to +\infty$}$$
where 
$C_{\beta,N,s}':= 2^{2s} \frac{ \Gamma\big(\frac{\beta}{2} +s\big)\Gamma\big(\frac{N}{2}-\frac{\beta}{2}\big)}{\Gamma\big(\frac{\beta}{2}\big)\Gamma\big(\frac{N}{2}-\frac{\beta}{2}-s\big) }<0.$

\item if $\beta=N-2s$, then 
\begin{align*}
(-\Delta)^s h_{N-2s}(x) &= C_{N-2s,N,s}' h_{N+2s}
(x) \quad \hbox{for $x \in \R^N\setminus \{0\}$}\\
& \sim C_{N-2s,N,s}' \frac{1}{|x|^{N + 2s}} \quad \hbox{as $|x|\to +\infty$}
\end{align*} 
where 
$C_{N-2s,N,s}':= 2^{2s} \frac{ \Gamma\big(\frac{N}{2} +s\big)}{\Gamma\big(\frac{N}{2}-s\big) } >0.$

\item if $\beta\in (0, N-2s)$, then 
$$(-\Delta)^s h_\beta(x) \sim C_{\beta,N,s}' \frac{1}{|x|^{\beta + 2s}} \quad \hbox{as $|x|\to +\infty$}$$
where 
$C_{\beta,N,s}':= 2^{2s} \frac{ \Gamma\big(\frac{\beta}{2} +s\big)\Gamma\big(\frac{N}{2}-\frac{\beta}{2}\big)}{\Gamma\big(\frac{\beta}{2}\big)\Gamma\big(\frac{N}{2}-\frac{\beta}{2}-s\big) }>0.$ 
This in particular includes the case $\beta = N - 2k$ with $k = 1, \dots, [\frac{N}{2}]$.
\end{itemize}
\end{Lemma}

\begin{Remark}
Notice that, for $\beta \in \{N-2s\} \cup [N, +\infty)$, the asymptotic behaviour of $\abs{(-\Delta)^s h_\beta(x)}$ does not depend on $\beta$; on the other hand, the sign and the precise constant depend on $\beta$.

In the case $\beta \in (0, N)\setminus \{N-2s\}$, we may use $x\mapsto \frac{1}{|x|^{\beta}}$, whose fractional Laplacian has a close (simple) representation:
$$\left((-\Delta)^s \frac{1}{|\cdot|^\beta} \right)(x) = C_{\beta,N,s} \frac{1}{|x|^{\beta+2s}},$$
see \cite[Table 1 and Theorem 3.1]{Kwa1}. 
In particular
$$(-\Delta)^s h_\beta(x) \sim \left((-\Delta)^s \frac{1}{|\cdot|^\beta} \right)(x) \quad \hbox{as $|x|\to +\infty$}.$$
On the other hand, if $\beta=N-2s$, we obtain, far from the origin, $(-\Delta)^s \frac{1}{|\,\cdot\,|^\beta} \equiv 0$ (recall that the Riesz potential $\frac{1}{|\cdot|^{N-2s}} \equiv I_{2s}$ is a fundamental solution, see e.g. \cite[Section 1.3.1]{GalT}); 
thus, in particular, the two functions have different asymptotic behaviours. 
This is the same reason why, for $h_{\beta}$, we have a discontinuity on the behaviour at infinity around $\beta=N-2s$.
\end{Remark}

\section*{Funding \& Acknowledgments}

The author is supported by PRIN 2017JPCAPN ``Qualitative and quantitative aspects of nonlinear PDEs'' and by INdAM-GNAMPA.

The author wishes to thank the Department of Mathematics of Università degli Studi di Bari Aldo Moro, where the paper was written (\cite{GalT}); in particular, he wishes to thank Prof. Silvia Cingolani for some fruitful discussions.


\begin{thebibliography}{999} 



\bibitem{AS0} M. Abramowitz, I. Stegun, 
``Handbook of Mathematical Functions with Formulas, Graphs, and Mathematical Tables'', 
Dover Publications Inc., 1964.

\bibitem{Amb0} V. Ambrosio,
\emph{Boundedness and decay of solutions for some fractional magnetic Schrödinger equations in $\R^N$},
Milan J. Math. \textbf{86} (2018), 125--136.

\bibitem{AAR} G.$\,$E. Andrews, R. Askey, R. Roy, 
``Special Functions'', 
Cambridge University Press, 1999.

\bibitem{ArMe} C. Argaez, M. Melgaard,
\emph{Solutions to quasi-relativistic multi-configurative Hartree–Fock equations in quantum chemistry},
Nonlinear Anal. \textbf{75} (2012), no. 1, 384--404.

\bibitem{BCI} G. Barles, E. Chasseigne, C. Imbert,
\emph{On the Dirichlet Problem for Second-Order Elliptic Integro-Differential Equations}, 
Indiana Univ. Math. J. \textbf{57} (2008), no. 1, 213--246.

\bibitem{BBMP} P. Belchior, H. Bueno, O.$\,$H. Miyagaki, G.$\,$A. Pereira,
\emph{Remarks about a fractional Choquard equation: ground state, regularity and polynomial decay}, 
Nonlinear Anal. \textbf{164} (2017), 38--53.

\bibitem{BL1} H. Berestycki, P.-L. Lions,
\emph{Nonlinear scalar field equations I. Existence of a ground state},
Arch. Rational Mech. Anal. \textbf{82} (1983), 313--345.


\bibitem{BBS} F. Bernini, B. Bieganowski, S. Secchi, 
\emph{Semirelativistic Choquard equations with singular potentials and general nonlinearities arising from Hartree–Fock theory},
Nonlinear Anal. \textbf{217} (2022), article ID 112738, pp. 26.


\bibitem{CafSil2} L.$\,$A. Caffarelli, L. Silvestre,
\emph{Regularity theory for fully nonlinear integro-differential equations}, 
Comm. Pure Appl. Math. \textbf{62} (2009), no. 5, 597--638.

\bibitem{CafSir} L.$\,$A. Caffarelli, Y. Sire, 
\emph{On some pointwise inequalities involving nonlocal operators}, 
in "Harmonic Analysis, Partial Differential Equations and Applications" (eds. S. Chanillo, B. Franchi, G. Lu, C. Perez, E.$\,$T. Sawyer), Springer, 2017.

\bibitem{CFQ} H. Chen, P. Felmer, A. Quaas,
\emph{Large solutions to elliptic equations involving fractional Laplacian},
Ann. Inst. H. Poincaré Anal. Non Linéaire C, \textbf{32} (2015), no. 6, 1199-1128.

\bibitem{CLO} W. Chen, C. Li, B. Ou, 
\emph{Classification of solutions for an integral equation},
Comm. Pure Appl. Math., \textbf{59} (2006), no. 3, 330--343.


\bibitem{CCS1} S. Cingolani, M. Clapp, S. Secchi,
\emph{Multiple solutions to a magnetic nonlinear Choquard equation},
Z. Angew. Math. Phys. \textbf{63} (2012), 233--248.

\bibitem{CG0} S. Cingolani, M. Gallo,
\emph{On the fractional NLS equation and the effects of the potential well's topology},
Adv. Nonlinear Stud. \textbf{21} (2021), no. 1, 1--40.

\bibitem{CG1} S. Cingolani, M. Gallo,
\emph{On some qualitative aspects for doubly nonlocal equations},
Discrete Contin. Dyn. Syst. Ser. S (2022), DOI: 10.3934/dcdss.2022041.

\bibitem{CGT2} S. Cingolani, M. Gallo, K. Tanaka,
\emph{Symmetric ground states for doubly nonlocal equations with mass constraint},
Symmetry \textbf{13} (2021), no. 7, article ID 1199, 1--17.

\bibitem{CGT3} S. Cingolani, M. Gallo, K. Tanaka,
\emph{On fractional Schr\"odinger equations with Hartree type nonlinearities}, 
Mathematics in Engineering \textbf{4} (2022), no. 6, 1–33.


\bibitem{CGT5} S. Cingolani, M. Gallo, K. Tanaka,
\emph{Infinitely many free or prescribed mass solutions for fractional Hartree equations and Pohozaev identities}, 
to appear on Adv. Nonlinear Stud. 
(2023), \href{https://doi.org/10.48550/arXiv.2305.14003}{arXiv:2305.14003}.

\bibitem{ClSa} M. Clapp, D. Salazar,
\emph{Positive and sign changing solutions to a nonlinear Choquard equation},
J. Math. Anal. Appl. \textbf{407} (2013), no. 1, 1--15.

\bibitem{DG0} L. D'Ambrosio, M. Ghergu,
\emph{Representation formulae for nonhomogeneous differential operators and applications to PDEs}, 
J. Differential Equations \textbf{317} (2022), 706--753.

\bibitem{DSS1} P. D'Avenia, G. Siciliano, M. Squassina,
\emph{On the fractional Choquard equations}, 
Math. Models Methods Appl. Sci. \textbf{25} (2015), no. 8, 1447--1476.



\bibitem{DOS} A. Dall'Acqua, T. {\O}stergaard S{\o}rensen, E. Stockmeyer,
\emph{Hartree-Fock theory for pseudorelativistic atoms}, 
Ann. Henri Poincar\'e \textbf{9} (2008), no. 4, 711--742.

\bibitem{DPV} E. Di Nezza, G. Palatucci, E. Valdinoci, 
\emph{Hitchhiker's guide to the fractional Sobolev spaces}, 
Bull. Sci. Math. \textbf{136} (2012), no. 5, 521--573.

\bibitem{DLWZPH} L. Dong, D. Liu, W. Qi, L. Wang, H. Zhou, P. Peng, C. Huang,
\emph{Necklace beams carrying fractional angular momentum in fractional systems with a saturable nonlinearity},
Commun. Nonlinear Sci. Numer. Simul. \textbf{99} (2021), article ID 105840, pp. 8.

\bibitem{DKK} B. Dyda, A. Kuznetsov, M. Kwaśnicki, 
\emph{Fractional Laplace operator and Meijer G-function},
Constr. Approx. \textbf{45} (2017), no. 3, 427--448.

\bibitem{ElSc0} A. Elgart, B. Schlein,
\emph{Mean field dynamics of boson stars},
Comm. Pure Appl. Math. \textbf{60} (2007), no. 4, 500--545.


\bibitem{FQT} P. Felmer, A. Quaas, J. Tan,
\emph{Positive solutions of the nonlinear Schr\"odinger equations with the fractional Laplacian},
Proc. Roy. Soc. Edinburgh Sect. A \textbf{142} (2012), no. 6, 1237--1262.

\bibitem{FLS} R.$\,$L. Frank, E. Lenzmann, L. Silvestre,
\emph{Uniqueness of radial solutions for the fractional Laplacian},
Comm. Pure Appl. Math. \textbf{69} (2016), no. 9, 1671--1726.

\bibitem{FJL} J. Fr\"ohlich, B.$\,$L.$\,$G. Jonsson, E. Lenzmann,
\emph{Boson stars as solitary waves}, 
Comm. Math. Phys. \textbf{274} (2007), no. 1, 1--30.
		

\bibitem{GalT} M. Gallo,
\emph{Nonlocal elliptic PDEs with general nonlinearities},
Ph.D. thesis, Università degli Studi di Bari Aldo Moro (2023).


\bibitem{Gar0} N. Garofalo, 
\emph{Fractional thoughts}, 
in "New Developments in the Analysis of Nonlocal Operators" \textbf{723} (eds. D. Danielli, A. Petrosyan, C.$\,$A. Pop), AMS, 2019.

\bibitem{Gre0} D. Greco, 
\emph{A Thomas-Fermi type variational problem with low regularity}, 
\href{https://arxiv.org/abs/2302.12586}{arXiv:2302.12586}, (2023).



\bibitem{GM0} G. Grillo, M. Muratori, 
\emph{On the asymptotic behaviour of solutions to the fractional porous medium equation with variable density},
Discrete Contin. Dyn. Syst. \textbf{35} (2015), no. 12, 5927--5962.


\bibitem{HLLS} C. Hainzl, E. Lenzmann, M. Lewin, B. Schlein,
\emph{On blowup for time-dependent generalized Hartree–Fock equations}, 
Ann. Henri Poincar\'e \textbf{11} (2010), no. 6, 1023--1052.


\bibitem{HL0} S. Herr, E. Lenzmann,
\emph{The Boson star equation with initial data of low regularity},
Nonlinear Anal. \textbf{97} (2014), 125--137.

\bibitem{JLX} T. Jin, Y. Y. Li, J. Xiong,
\emph{On a fractional Nirenberg problem, part I: blow up analysis and compactness of solutions},
J. Eur. Math. Soc. \textbf{16}, (2014), 1111–1171.

	
\bibitem{Kwa1} M. Kwaśnicki, 
\emph{Fractional Laplace operator and its properties}, 
in: A. Kochubei, Y. Luchko, "Handbook of Fractional Calculus with Applications, Vol. 1: Basic Theory", De Gruyter, Berlin, 2019.


\bibitem{Le0} P. Le,
\emph{Liouville theorem and classification of positive solutions for a fractional Choquard type equation}, 
Nonlinear Anal. \textbf{185} (2019), 123--141. 

\bibitem{Le1} P. Le,
\emph{Symmetry of positive solutions to Choquard type equations involving the fractional $p$-Laplacian}, 
Acta Appl. Math. \textbf{170} (2020), no. 1, 387--398. 



\bibitem{Lie2} E.$\,$H. Lieb,
\emph{Sharp constants in the Hardy-Littlewood-Sobolev and related inequalities}, 
Annals of Math. \textbf{118} (1983), no. 2, 349--374. 
	


\bibitem{LMM} J. Lu, V. Moroz, C.$\,$B. Muratov,
\emph{Orbital-free density functional theory of out-of-plane charge screening in graphene},
J. Nonlinear. Sci. \textbf{25} (2015), no. 6, 1391--1430.

\bibitem{MZ0} P. Ma, J. Zhang,
\emph{Existence and multiplicity of solutions for fractional Choquard equations},
Nonlinear Anal. \textbf{164} (2017), pp. 100--117.


\bibitem{MPS2} L. Maia, B. Pellacci, D. Schiera,
\emph{Symmetric positive solutions to nonlinear Choquard equations with potentials},
Calc. Var. Partial Differential Equations \textbf{61} (2022), no. 61, pp. 34.

\bibitem{MS0} V. Moroz, J. Van Schaftingen, 
\emph{Groundstates of nonlinear Choquard equations: existence, qualitative properties and decay asymptotics}, 
J. Funct. Anal. \textbf{265} (2013), no. 2, 153--184.


\bibitem{MS2} V. Moroz, J. Van Schaftingen, 
\emph{Existence of groundstates for a class of nonlinear Choquard equations}, 
Trans. Amer. Math. Soc., \textbf{367} (2015), no. 9, 6557--6579. 

\bibitem{MS3} V. Moroz, J. Van Schaftingen, 
 \emph{A guide to the Choquard equation}, 
 J. Fixed Point Theory Appl. \textbf{19} (2017), no. 1, 773--813.

	

\bibitem{PT0} F. Punzo, G. Terrone, 
\emph{On a fractional sublinear elliptic equation with a variable coefficient},
Appl. Anal. \textbf{94} (2015), no. 4, 800--818.

\bibitem{QX0} A. Quaas, A. Xia, 
\emph{Liouville type theorems for nonlinear elliptic equations and systems involving fractional Laplacian in the half space},
Calc. Var. Partial Differential Equations \textbf{52} (2015), 641--659.

\bibitem{ROS1} X. Ros-Oton, J. Serra,
\emph{The Dirichlet problem for the fractional Laplacian: regularity up to the boundary},
J. Math. Pures Appl. \textbf{101} (2014), no. 3, 275--302.

\bibitem{SeV} R. Servadei, E. Valdinoci,
\emph{Weak and viscosity solutions of the fractional Laplace equation},
Publ. Mat. \textbf{58} (2015), no. 1, 133--154.

\bibitem{ShSp} T.-T. Shieh, D.$\,$E. Spector,
\emph{On a new class of fractional partial differential equations},
Adv. Calc. Var. \textbf{8} (2014), no. 4, pp. 16.

\bibitem{Sil0} L. Silvestre,
\emph{Regularity of the obstacle problem for a fractional power of the Laplace operator},
Comm. Pure Appl. Math. \textbf{60} (2006), 67--112.

\bibitem{SoV} N. Soave, E. Valdinoci,
\emph{Overdetermined problems for the fractional Laplacian in exterior and annular sets},
J. Anal. Math. \textbf{137} (2019), no. 1, 101--134.



\bibitem{Tri0} H. Triebel,
``Interpolation theory, function spaces, differential operators",
North-Holland \textbf{18}, 1978.


\bibitem{WY0} X. Wang, Z. Yang,
\emph{Symmetry and monotonicity of positive solutions for a Choquard equation with the fractional Laplacian}, 
Complex Var. Elliptic Equ. \textbf{67} (2021), no. 5, 1211-1228.	

\bibitem{WG0} Z.$\,$X. Wang, D.$\,$R. Guo,
``Special Functions'', 
World Scientific, Singapore, 2010.

\bibitem{ZY0} S. Zhao, Y. Yu,
\emph{Sign-changing solutions for a fractional Choquard equation with power nonlinearity}, 
Nonlinear Anal. \textbf{221} (2022), article ID 112917, pp. 18.

\end{thebibliography}
\end{document}